\documentclass{amsart}
\usepackage{amsmath}
\usepackage{amssymb} 
\usepackage{amscd} 
\usepackage{graphicx} 
\usepackage{epsfig}
\usepackage[dvips]{color}
\usepackage{bm}

\newtheorem{theorem}{Theorem}[section]
\newtheorem{lemma}[theorem]{Lemma}

\newtheorem{proposition}[theorem]{Proposition}
\newtheorem{corollary}[theorem]{Corollary}
\newtheorem{claim}[theorem]{Claim}

\theoremstyle{definition}
\newtheorem{definition}[theorem]{Definition}

\theoremstyle{remark}
\newtheorem{remark}[theorem]{Remark}

\newtheorem{question}[theorem]{Question}

\numberwithin{equation}{section}
\numberwithin{figure}{section}
\numberwithin{table}{section}

\newcommand{\QED}[1]{\hspace*{\fill} $\square$(#1)\newline}

\begin{document}

\title[Non primitive/Seifert-fibered Dehn surgeries]
{Seifert fibered surgeries on strongly invertible knots without primitive/Seifert positions}

\author[M. Eudave-Mu\~noz]{Mario Eudave-Mu\~noz}
\address{Instituto de Matematicas, Universidad Nacional Aut\'onoma de M\'exico, Circuito Exterior, Ciudad Universitaria 04510 M\'exico DF, Mexico, and CIMAT, Guanajuato, Mexico}
\email{mario@matem.unam.mx}
\thanks{The first author was partially supported by PAPIIT-UNAM grant IN109811.}

\author[E. Jasso]{Edgar Jasso}
\address{Mathematics Department, 
North Seattle Community College, 
9600 College Way N, 
Seattle, WA 98103, USA}
\email{ejasso@northseattle.edu}

\author[K. Miyazaki]{Katura Miyazaki}
\address{Faculty of Engineering, Tokyo Denki University, 5 Senju Asahi-cho, Adachi-ku, Tokyo 120-8551, 
Japan}
\email{miyazaki@cck.dendai.ac.jp}

\author[K. Motegi]{Kimihiko Motegi}
\address{Department of Mathematics, Nihon University, 
3-25-40 Sakurajosui, Setagaya-ku, 
Tokyo 156-8550, Japan}
\email{motegi@math.chs.nihon-u.ac.jp}
\thanks{
The last author has been partially supported by JSPS Grants-in-Aid for Scientific 
Research (C) (No.21540098), The Ministry of Education, Culture, Sports, Science and Technology, Japan and Joint Research Grant of Institute of Natural Sciences at 
Nihon University for 2013. }

\subjclass{Primary 57M25, 57M50 Secondary 57N10}
\date{}

\keywords{Dehn surgery, hyperbolic knot, Seifert fiber space, primitive/Seifert position, 
seiferter, Seifert Surgery Network, branched covering, Montesinos trick}

\dedicatory{Dedicated to Michel Boileau on the occasion of his 60th birthday}

\begin{abstract}
We find an infinite family of Seifert fibered surgeries on strongly invertible knots which do not have primitive/Seifert positions. 
Each member of the family is obtained from a trefoil knot after alternate twists along a pair of seiferters for 
a Seifert fibered surgery on a trefoil knot.  
\end{abstract}

\maketitle

\section{Introduction}
\label{section:Introduction}
A pair $(K, m)$ of a knot $K$ in the $3$--sphere $S^3$ and an integer $m$ is called 
a \textit{Seifert fibered surgery} if the resulting manifold $K(m)$ obtained by $m$--surgery on $K$ is a Seifert fiber space. 
For most known Seifert fibered surgeries $(K, m)$, 
$K$ has a nice position on the boundary of
a standard genus 2 handlebody in $S^3$,
which is called a ``primitive/Seifert position".

For a genus $2$ handlebody $H$ and a simple closed curve $c$ in
$\partial H$, we denote $H$ with a $2$--handle attached along $c$
by $H[c]$.
Let $S^3 = V \cup_F W$ be a genus $2$ Heegaard splitting of $S^3$, 
i.e.\ $V$ and $W$ are genus $2$ handlebodies in $S^3$
with $V \cap W$ a genus $2$ Heegaard surface $F$. 
We say that a Seifert fibered surgery $(K, m)$ has a \textit{primitive/Seifert position} if 
there is a genus $2$ Heegaard surface $F$ which carries $K$ and satisfies the following three conditions. 

\begin{itemize}
\item 
$K$ is primitive with respect to $V$, 
i.e.\ $V[K]$ is a solid torus. 
\item
$K$ is Seifert with respect to $W$, 
i.e.\ $W[K]$ is a Seifert fiber space 
over the disk with two exceptional fibers.
\item
The surface slope of $K$ with respect to $F$ 
(i.e.\ the isotopy class in  $\partial N(K)$ represented by a component of $\partial N(K) \cap F$) 
coincides with the surgery slope $m$. 
\end{itemize}

Assume that a knot $K$ has a primitive/Seifert position
with surface slope $m$. 
Then $K(m) \cong V[K] \cup W[K]$ is
a Seifert fiber space or a connected sum of lens spaces. 
Moreover, $K$ has tunnel number one \cite[2.3]{D}, 
and hence strongly invertible \cite[Claim 5.3]{MM7}. 
By the positive solution to the cabling conjecture for strongly
invertible knots \cite{EM},
if $K$ is hyperbolic,
then $K(m)$ is a Seifert fiber space over $S^2$ with
at most three exceptional fibers,
so $(K, m)$ is a Seifert fibered surgery.

Primitive/Seifert positions,
introduced by Dean \cite{D}, are variants 
of Berge's primitive/primitive positions \cite{Berge2}. 
Although any lens surgery is conjectured to have
a primitive/primitive position \cite{Berge2, Go90}, 
there are infinitely many Seifert fibered surgeries 
with no primitive/Seifert positions \cite{MMM, DMM1, Tera}. 
Knots yielding these Seifert fibered surgeries
are not strongly invertible;
the simplest example is $1$--surgery on the pretzel knot $P(-3, 3, 5)$. 
So it is natural to ask:

\begin{question}
\label{question1}
Let $(K, m)$ be a Seifert fibered surgery 
on a strongly invertible knot $K$.
Then, does it have a primitive/Seifert position?
\end{question}

However,
Song \cite{Song} observed that $1$--surgery on $P(-3, 3, 3)$
yields a Seifert fiber space.
Since $P(-3, 3, 3)$ is a strongly invertible
knot of tunnel number $2$,
that surgery gives the negative answer to Question~\ref{question1}.
In \cite{DMMtrefoil} we construct a one--parameter family
of Seifert fibered surgeries which
answer Question~\ref{question1} in the negative
and contain Song's example
by using the Seifert Surgery Network
introduced in \cite{DMM1}.
In this paper,
we construct a large family of Seifert fibered surgeries
giving the negative answer to Question~\ref{question1}
by taking $2$--fold branched covers of tangles. 
We then study these surgeries
from a viewpoint of the Seifert Surgery Network,
and find a path in the network from each surgery in our family to
a surgery on a trefoil knot. 
Our family of Seifert fibered surgeries
is a variant of families obtained in \cite{EM2, DEMM}.
In \cite{EM2}, 
by using $2$--fold branched covers of tangles,
the first author constructs 4 families of Seifert fibered surgeries
having primitive/Seifert positions. \par

We briefly review the definitions of seiferters and
the Seifert Surgery Network.
For a knot $K \subset S^3$ and $m \in \mathbb{Z}$,
the pair $(K, m)$ is a \textit{Seifert surgery} if
$K(m)$ has a possibly degenerate Seifert fibration,
i.e.\ a Seifert fibration which may contain
an exceptional fiber of index $0$.
Let $(K, m)$ be a Seifert surgery.
A simple closed curve $c$ in $S^3 -K$ is called a \textit{seiferter}
if $c$ is a trivial knot in $S^3$ and
a Seifert fiber in $K(m)$.
Denoting by $K_p$ and $m_p$
the images of $K$ and $m$ under $p$--twist along the seiferter $c$,
we see that $(K_p, m_p)$ is a Seifert surgery
with $c$ a seiferter.
If seiferters $c_1, c_2$ for $(K, m)$
become fibers in a Seifert fibration of $K(m)$ simultaneously,
then $\{c_1, c_2\}$ is a \textit{pair of seiferters} for $(K, m)$.
If a pair of seiferters cobound an annulus $A$ in $S^3$,
the pair is an \textit{annular pair of seiferters}.
As is twisting along a seiferter,
twisting $(K, m)$ along the annulus $A$ yields a Seifert surgery. 
\textit{The Seifert Surgery Network} is
the $1$--dimensional complex such that
its vertices are Seifert surgeries and
two vertices are connected by an edge if
one is obtained from the other by $1$--twist along a seiferter
or an annular pair of seiferters;
see \cite[Subsection~2.4]{DMM1}. 
Hence, a path in the Seifert Surgery Network tells
how one Seifert surgery is obtained from another 
by twisting along seiferters and/or
annular pairs of seiferters.

Our main result is as follows.

\begin{figure}[h]
\begin{center}
\includegraphics[width=0.3\linewidth]{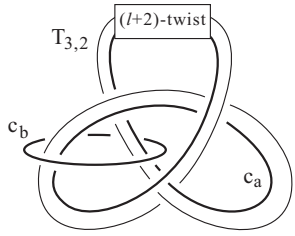}
\caption{$\{c_a, c_b\}$ is a pair of seiferters for $(T_{3,2}, l+5)$.}
\label{fig:c_a c_b}
\end{center}
\end{figure}

\begin{theorem}
\label{non p/s}
There are infinitely many Seifert fibered surgeries on strongly invertible
hyperbolic knots 
$(K(l, m, n, p), \gamma_{l,m,n,p})$ 
$(m=0$ or $p=0)$
with the following properties,
where $l,m,n,p$ satisfy more conditions given
in Proposition~\ref{nonPS}. 
\begin{enumerate}
\item 
$(K(l, m, n, p), \gamma_{l,m,n,p})$ does not have a primitive/Seifert position.
\item
The pair of knots $c_a, c_b$ in Figure~\ref{fig:c_a c_b}
is a pair of seiferters for $(l + 5)$--surgery
on the trefoil knot $T_{3, 2}$.
$(K(l, m, n, p), \gamma_{l,m,n,p})$ is obtained from
$(T_{3, 2}, l + 5)$ by applying a sequence of twists
along $c_a, c_b$.
Refer to Proposition~\ref{location of Klmnp} and also
Corollary~\ref{cor:location of Klmnp} for the details of the sequence. 

\end{enumerate}
\end{theorem}

\textsc{Proof of Theorem~\ref{non p/s}.}
The definition of $K(l, m, n, p)$ and $\gamma_{l,m,n,p}$
is given in Subsection~\ref{subsection:Klmnp}.
Proposition~\ref{Klmnp} shows that
$K(l,m,n,p)(\gamma_{l,m,n,p})$ is a Seifert fiber space.
The hyperbolicity of $K(l,m,n,p)$ is proved in
Proposition~\ref{hyperbolic}.
Assertion~(1) follows from Proposition~\ref{nonPS}
by showing that the tunnel number of $K(l,m,n,p)$ is 2.
Assertion~(2) follows from Proposition~\ref{location of Klmnp}.
\QED{Theorem~\ref{non p/s}}

For all the known Seifert fibered surgeries $(K, m)$ with
no primitive/Seifert positions, 
$K$ has tunnel number greater than one.
We close with the following question. 

\begin{question}
\label{question2}
Let $(K, m)$ be a Seifert fibered surgery
 on a tunnel number one knot $K$. 
Then, does it have a primitive/Seifert position?
\end{question}

\section{Tangles, branched coverings and Seifert fibered surgeries}
\label{section:covering}

Let $B$ be a $3$--ball and $t$ a disjoint union of two arcs properly embedded in $B$ and some 
simple closed curves. 
Then the pair $(B, t)$ is called a \textit{tangle}. 
A tangle $(B, t)$ is trivial if there is a pairwise
homeomorphism from $(B, t)$ to $(D^2 \times I, \{x_1, x_2\}\times I)$, where $x_1, x_2$ are distinct points.

\begin{figure}[h]
\begin{center}
\includegraphics[width=0.9\linewidth]{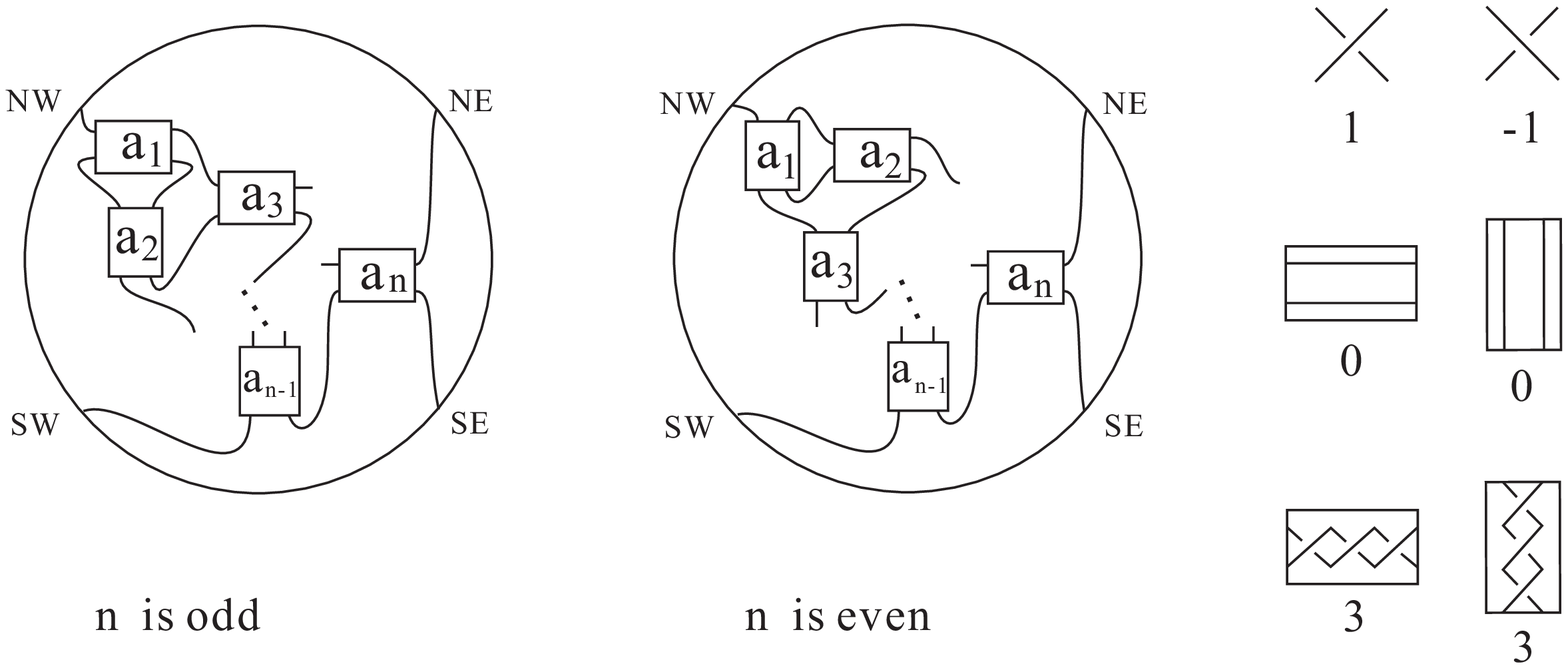}
\caption{Rational tangles.}
\label{rtangle}
\end{center}
\end{figure}

Let $U$ be the unit 3-ball in $\mathbb{R}^3$,
and take 4 points NW, NE, SE, SW on the boundary of $U$
so that
$\mathrm{NW} = (0,-\alpha, \alpha),
\mathrm{NE} = (0, \alpha, \alpha),
\mathrm{SE} = (0, \alpha, -\alpha),
\mathrm{SW} = (0, -\alpha, -\alpha)$,
where $\alpha =\frac{1}{\sqrt{2}}$.
A tangle $(U, t)$ is a \textit{rational tangle}
if it is a trivial tangle with $\partial t =
\{\mathrm{NW, NE, SE, SW}\}$. 
Two rational tangles $(U, t)$ and $(U, t')$ are \textit{equivalent} if there is 
a pairwise homeomorphism $h : (U, t) \to (U, t')$ such that $h|_{\partial U}$ is the identity map. 
We can construct rational tangles from sequences of integers 
$a_1, a_2, \dots, a_n$ as shown in Figure~\ref{rtangle}, 
where the last horizontal twist $a_n$ may be $0$. 
In our figure a horizontal rectangle (i.e. a rectangle intersecting arcs $t$ on the left and right sides) 
with a label $a_i$ represents a strand of $a_i$ horizontal crossings, 
with the sign convention shown in Figure~\ref{rtangle}. 
Similarly, 
a vertical rectangle  (i.e. a rectangle intersecting arcs $t$ on the top and bottom sides) with a label $a_i$ represents a strand of $a_i$ vertical crossings, with the sign convention shown in Figure~\ref{rtangle}. 
We consider that the tangle diagrams in Figure~\ref{rtangle}
are drawn on the $yz$--plane. 
Denote by $R(a_1, a_2, \dots, a_n)$ the associated rational tangle.  

Each rational tangle can be parametrized by 
$r \in \mathbb{Q} \cup \{ \infty\}$, 
where the rational number $r$ is given by
the continued fraction below. 
Thus we denote the rational tangle corresponding to $r$ by $R(r)$. 

$$ r\ =\ a_{n} + \cfrac 1 {a_{n-1}+ \cfrac 1 {
               \begin{array}{clr}
               \ & & \\[-5pt]
               \hspace{-25pt} \ddots & & \\[-10pt]
                      & \ \ \hspace{-5pt} a_2+ \cfrac{1}{a_1}
               \end{array}
               }}$$
               
\begin{figure}[h]
\begin{center}
\includegraphics[width=0.35\linewidth]{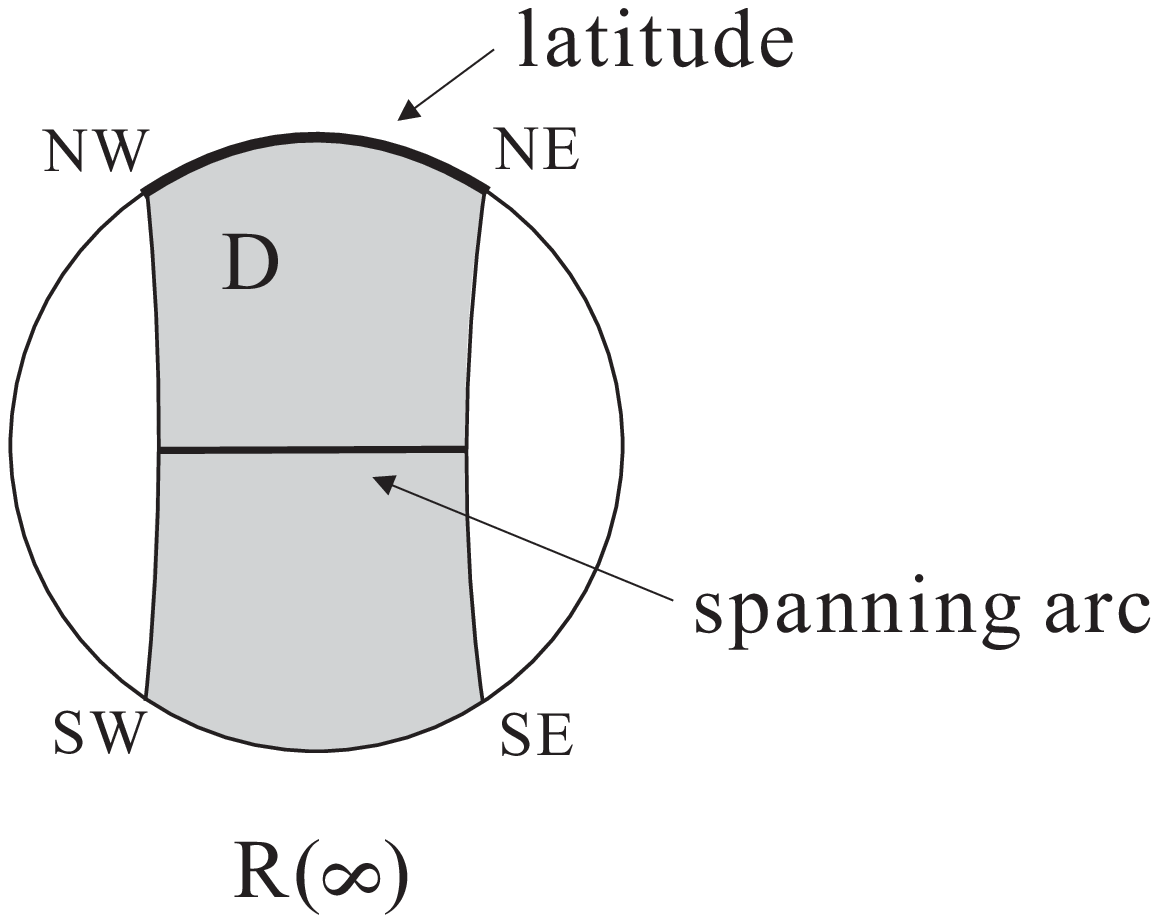}
\caption{}
\label{spanningarc}
\end{center}
\end{figure}

Let $(U, t)$ be the rational tangle $R(\infty)$.
Considering $t$ is embedded in the $yz$--plane,
take the disk $D$ in the $yz$--plane such that
$\partial D$ is the union of $t$ and 
two arcs in $\partial U$: one connects
NW and NE, and the other connects SW and SE.
We call an arc in $D$ connecting the components of the interior of $t$
a \textit{spanning arc},
and the arc $D\cap \partial U$ connecting
NW and NE \textit{the latitude of $R(\infty)$}. 
See Figure~\ref{spanningarc}. 
The $2$--fold cover $\widetilde{U}$
of $U$ branched along $t$ is a solid torus.
Note that the preimages of the spanning arc
and the latitude are the core and a longitude $\lambda$
of the solid torus, respectively.
A \textit{meridian of a rational tangle} $R(r)=(U, t')$ is
a simple closed curve in $\partial U -t'$ which bounds a disk
in $U -t'$ and a disk in $\partial U$ meeting $t'$
in two points.
Let $\mu_r (\subset \partial \widetilde{U})$
be a lift of a meridian of $R(r)$;
then $\mu_r$ is a meridian of the solid torus $\widetilde{U}$.
Furthermore, we note the following well-known fact.

\begin{lemma}
\label{-p/q}
Under adequate orientations we have
$[\mu_r] = -p[\mu_{\infty}] +q[\lambda]
\in H_1(\partial \widetilde{U})$,
where $r = \frac{p}{q}$ and $[\mu_{\infty}]\cdot[\lambda] =1$.
\end{lemma}

Let $(B, t)$ be a tangle such that
$B \subset S^3 ( = \mathbb{R}^3\cup\{\infty\} )$ is
the complement of the unit 3-ball $U$,
and $\partial t = \{ \mathrm{NW, NE, SE, SW} \}$.
We denote by $(B, t) + R(r)$ 
the knot or link in $S^3$ 
formed by the union of the strings of the tangles,
and let $\pi_r: X_r \to S^3 =B\cup U$ be
the $2$--fold cover branched along $(B, t) + R(r)$.
We say that $(B, t)$ is \textit{trivializable}
if  $(B, t) + R(\infty)$ is a trivial knot in $S^3$. 
If $(B, t) + R(r)$ is a trivial knot for some $r\in \mathbb{Q}$, 
then an ambient isotopy of $B$ changes
$(B, t)$ to a trivializable tangle.

Suppose that $(B, t)$ is trivializable.
Then the $2$--fold branched cover $X_{\infty}$ is the $3$--sphere,
and the preimage of the spanning arc $\kappa$ for $R(\infty)$
is a knot in $X_{\infty} =S^3$,
which we call the \textit{covering knot} of $(B, t)$.
Note that the covering knot is a strongly invertible knot
whose strong inversion is the covering transformation of $X_{\infty}$.
The exterior of the covering knot $K$
is $\pi_{\infty}^{-1}(B)$.
For $(B, t) +R(\infty)$
a replacement of $R(\infty)$ by a rational tangle $R(s)$
is called \textit{$s$--untangle surgery}
on $(B, t) +R(\infty)$. 
Performing untangle surgery downstairs
corresponds to replacing the solid torus $\pi_{\infty}^{-1}(U)$
by $\pi_s^{-1}(U)$ upstairs,
i.e.\ Dehn surgery on the covering knot $K$.
We denote the surgery slope by $\gamma_s$;
it is represented by a lift
of a meridian of $R(s)$. 
We say that $\gamma_s$ is the \textit{covering slope} of $s$. 
See the commutative diagram below. 

\begin{eqnarray*}
\begin{CD}
	S^3 @>\gamma_s\textrm{--surgery on } K >> K(\gamma) \\
@V{ \textrm{2--fold branched cover}}VV
		@VV{\textrm{2--fold branched cover}}V \\
	(B, t)\cup R(\infty) @>>s\textrm{--untangle surgery} > (B, t)\cup R(s)
\end{CD}
\end{eqnarray*}

\begin{center}
\textsc{Diagram 2.} Montesinos trick 
\end{center}

\begin{remark}
\label{n--framing}
Suppose that the preimage of the latitude of $R(\infty)$
is a longitude of $\pi_{\infty}^{-1}( U )$ giving an $n$--framing.
Then, by Lemma~\ref{-p/q}
the covering slope $\gamma_s$ is $n-s$ 
in terms of a preferred meridian--longitude pair of $K$.
\end{remark}

A \textit{sum of two tangles} $(B_1, t_1)$ and $(B_2, t_2)$ is
the knot or link obtained by attaching $t_1$ and $t_2$ via
an orientation reversing homeomorphism
$h: \partial B_1 \to \partial B_2$
with $h(\partial t_1) =\partial t_2$.

\begin{figure}[h]
\begin{center}
\includegraphics[width=0.75\linewidth]{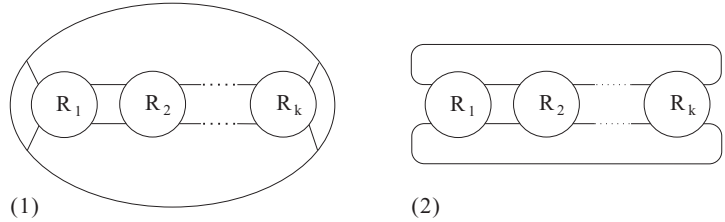}
\caption{Montesinos tangle and Montesinos link.}
\label{Montesinos}
\end{center}
\end{figure}

For rational tangles $R_1, \dots, R_k$,
the tangle in Figure~\ref{Montesinos}(1) is called 
a \textit{Montesinos tangle} $M_{T}(R_1, \ldots, R_k)$. 
The knot or link in Figure~\ref{Montesinos}(2) is
called a \textit{Montesinos link} $M(R_1, \ldots, R_k)$. 
We call the diagrams in Figure~\ref{Montesinos}
\textit{standard positions} of a Montesinos tangle
and a Montesinos link. 
If $R_i$ corresponds to $r_i\in \mathbb{Q}\cup \{\infty\}$ for $i=1, \ldots, k$,
then we often write
$M_T(r_1, \ldots, r_k)$ and $M(r_1, \ldots r_k)$ 
for a Montesinos tangle and a Montesinos link, respectively.

\begin{figure}[h]
\begin{center}
\includegraphics[width=0.65\linewidth]{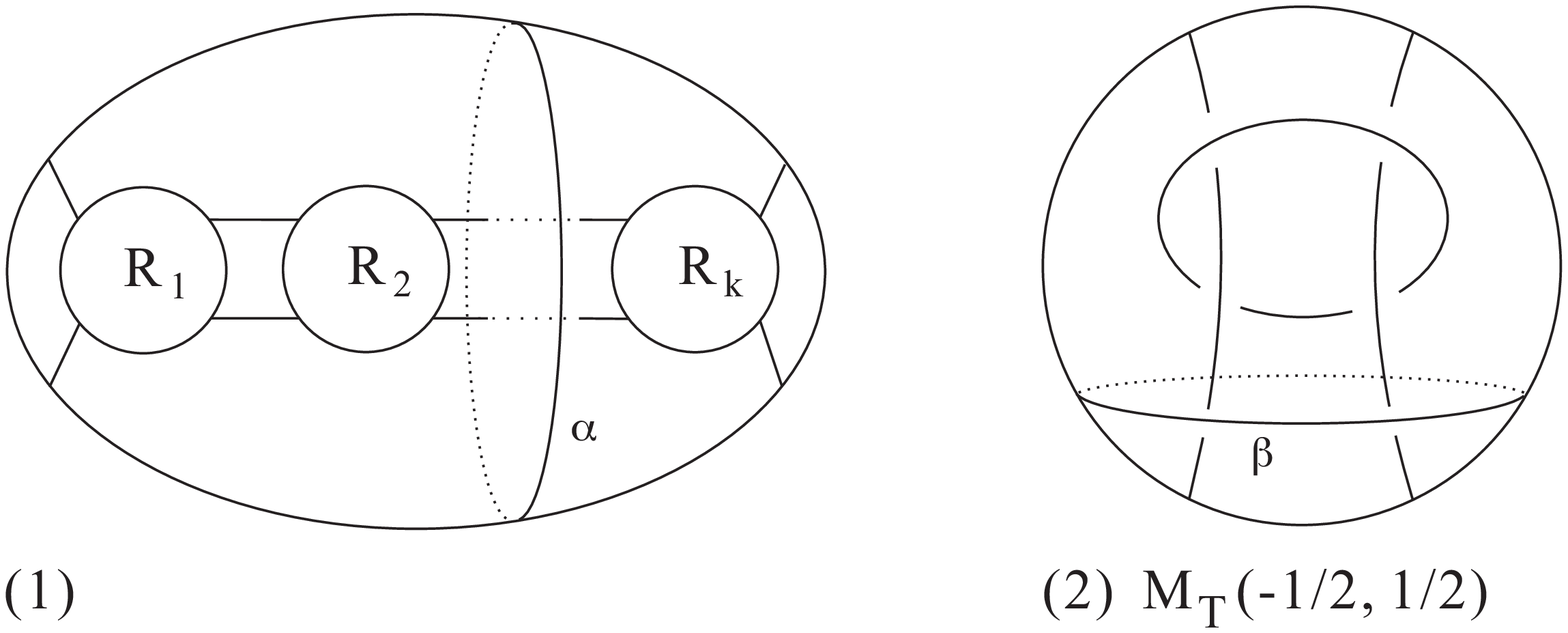}
\caption{}
\label{Mtangle_fiber}
\end{center}
\end{figure}

Let $X$ be the $2$--fold branched cover of $S^3$ (resp.\ $D^3$) along
a Montesinos link $M(R_1, \ldots, R_k)$ (resp.\ 
a Montesinos tangle $M_T(R_1, \ldots, R_k)$).
Then $X$ admits a Seifert fibration over $S^2$ (resp.\ $D^2$) in which
the preimage of $B_i$, 
where $R_i =(B_i, \frac{p_i}{q_i})$,
is a fibered solid torus and its core is an exceptional fiber
of Seifert invariant $\frac{p_i}{q_i}$ and index $|q_i|$.
Hence, $X$ is a Seifert fiber space $S^2(r_1, \ldots, r_k)$ 
(resp. $D^2(r_1, \ldots, r_k)$),
where $r_i = \frac{p_i}{q_i}$.  
Note that for the $2$--fold branched cover of $D^3$
along $M_T(r_1, \ldots, r_k)$,
a lift of a simple closed curve $\alpha (\subset \partial D^3)$
in Figure~\ref{Mtangle_fiber}(1)
is a fiber of $X =D^2(r_1, \ldots, r_k)$ up to isotopy.
See \cite{Montesinos}.

We have Lemma~\ref{lemma:fiber} below on
Seifert fibrations of $D^2(r_1, \ldots, r_k)$.

\begin{lemma}
\label{lemma:fiber}
Let $X$ be the $2$--fold branched cover of the $3$--ball along
a Montesinos tangle $M_T(\frac{p_1}{q_1}, \ldots, \frac{p_k}{q_k})$,
where $|q_i| \ge 2$ for all $i$.
Then the following hold.
\begin{enumerate}
\item $X$ admits more than one Seifert fibrations up to isotopy
if and only if $k = 2$ and $|q_1| = |q_2| = 2$.

\item Assume $k=2$ and $|q_1| =|q_2| =2$.
Then, $X$ is the twisted $S^1$ bundle over the M\"obius band,
and admits exactly two Seifert fibrations:
one is over the disk, and the other is over the M\"obius band
with no exceptional fibers.
In the latter fibration,
a fiber on $\partial X$ is isotopic in $\partial X$ to
a lift of a simple closed curve $\beta$ on $\partial B$
in Figure~\ref{Mtangle_fiber}(2),
where $(B, t) = M_T(-\frac{1}{2}, \frac{1}{2})$.
\end{enumerate}
\end{lemma}

\noindent
\textsc{Proof of Lemma~\ref{lemma:fiber}.}
We only prove the last statement of (2).
Let $\pi: X \to B$ be the $2$--fold cover branched along $t$,
where $(B, t) = M_T(-\frac{1}{2}, \frac{1}{2})$.
Let $A$ be an annulus properly embedded in $B -t$
such that a component of $\partial A$ is $\beta$
and $A$ separates the circle component from the two arcs of $t$. 
Then $\pi^{-1}(A)$ consists of two annuli and
splits $X$ into two solid tori.
Each component of $\pi^{-1}(A)$ is a non-separating
annulus in $X$, and thus a vertical annulus (i.e.\ a union of fibers)
in a Seifert fibration of $X$ over the M\"obius band.
This implies the claimed result.
\QED{Lemma~\ref{lemma:fiber}}

Let $(B, t)$ be a trivializable tangle such that
$(B, t) + R(s)$ is a Montesinos link
for some rational number $s$.
The $2$--fold branched cover $X_s$,
which is a Seifert fiber space as shown above,
is obtained from $S^3$ by
$\gamma_s$--surgery on the covering knot $K$ of $(B, t)$.
In this manner, we obtain a Seifert fibered surgery $(K, \gamma_s)$.

\section{Non primitive/Seifert-fibered, Seifert fibered surgeries on covering knots}
\label{section:examples}

In Subsection~\ref{subsection:Klmnp},
we construct Seifert surgeries $(K(l, m, n, p), \gamma_{l,m,n,p})$ ($m=0$ or $p = 0$) 
by untangle surgeries of trivializable tangles.
The knots $K(l, m, n, p)$ are strongly invertible.
In Subsection~\ref{Non primitive/Seifert},
these surgeries are shown to have no primitive/Seifert positions. 
In Subsection~\ref{subsection:hyperbolicity},
$K(l,m,n,p)$ are shown to be hyperbolic knots.

\subsection{Seifert fibered surgeries on knots \bm{$K(l,m,n,p)$}}
\label{subsection:Klmnp}

Let $B(l, m, n, p)$ be the tangle of Figure~\ref{Blmnp}. 
Then we have Lemma~\ref{Blmnp+R} below.

\begin{figure}[h]
\begin{center}
\includegraphics[width=0.55\linewidth]{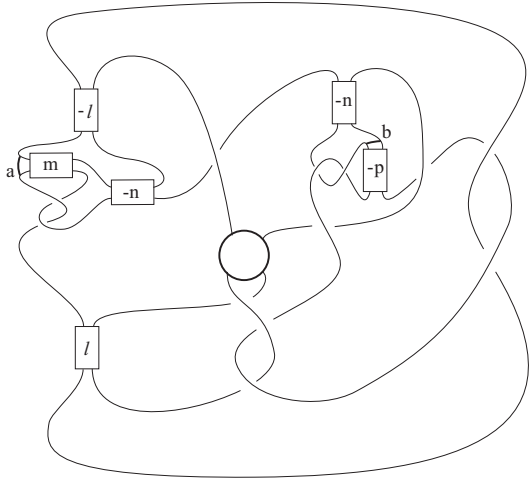}
\caption{Tangle $B(l, m, n, p)$: $m$ or $p$ is zero.}
\label{Blmnp}
\end{center}
\end{figure}

\begin{lemma}
\label{Blmnp+R}
\begin{enumerate}
\item
The tangle $B(l, m, n, 0)$ enjoys the following properties. 
	\begin{enumerate}
	\item $B(l, m, n, 0) + R(\infty)$ is a trivial knot.
	
	\smallskip	
	\item $B(l, m, n, 0) + R(1)$ is the Montesinos link 
	
	$\displaystyle
	M(\frac{2lmn+lm-ln+2mn+3m-n-1}{2l^2mn+l^2m-l^2n+2lm-2m-l+1}, 
-\frac{n+1}{4n+3}, \frac{1}{2})$. 
	\end{enumerate}
	
\smallskip
\item
The tangle $B(l, 0, n, p)$ enjoys the following properties. 
	\begin{enumerate}
	\item
	$B(l, 0, n, p) + R(\infty)$ is a trivial knot.
	
	\smallskip
	\item 
	$B(l, 0, n, p) + R(1)$ is the Montesinos link 
	\smallskip
	
	$\displaystyle
	M(\frac{ln+n+1}{l^2n+l-1}, \frac{-2np+n-p+1}{8np-4n+2p-3},
\frac{1}{2})$. 
	
	\end{enumerate}
\end{enumerate}
\end{lemma}

\begin{figure}[h]
\begin{center}
\includegraphics[width=1.0\linewidth]{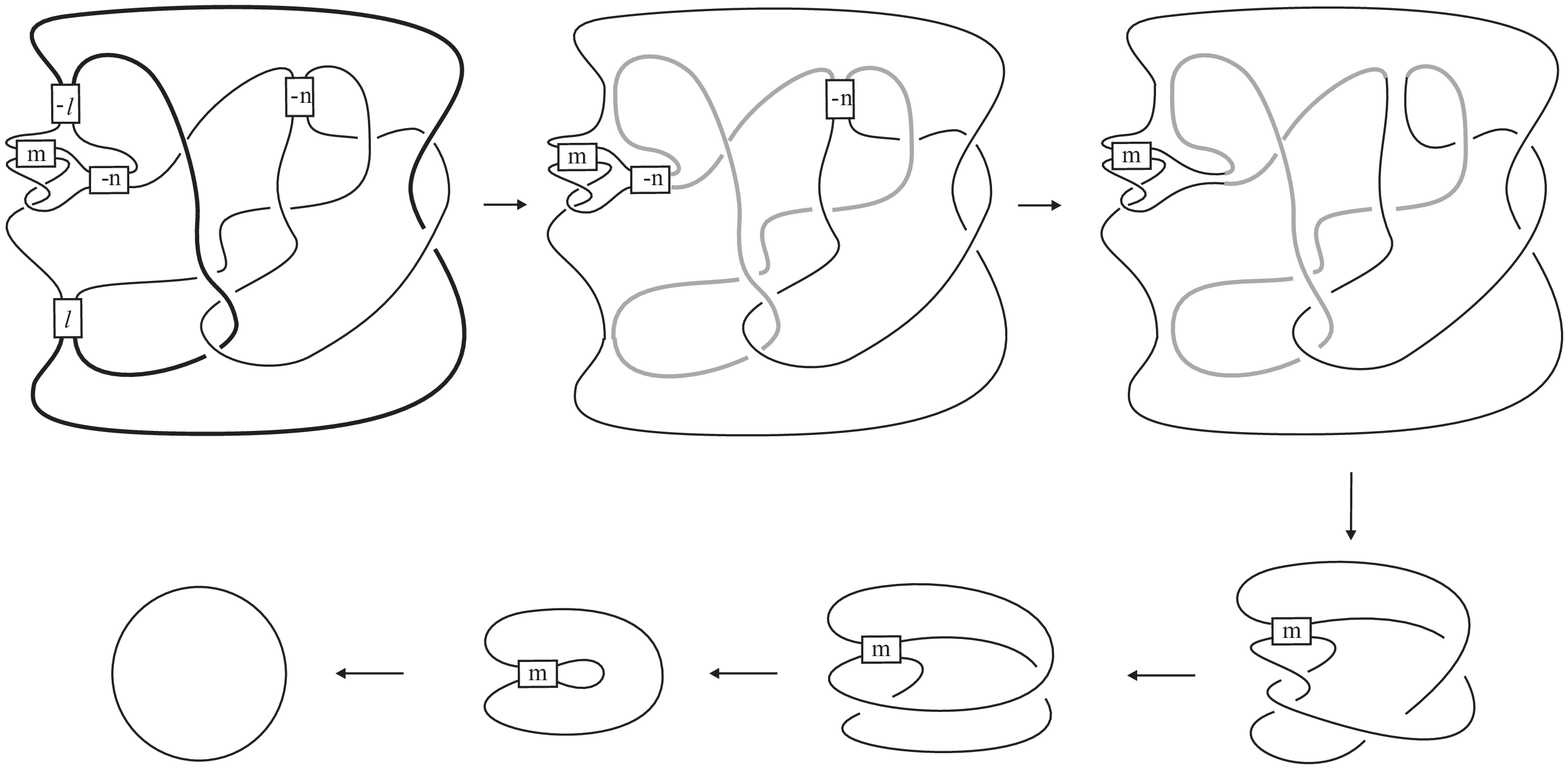}
\caption{ $B(l, m, n, 0) + R(\infty)$ is a trivial knot.}
\label{Blmn0+8unknot}
\end{center}
\end{figure}

\begin{figure}[h]
\begin{center}
\includegraphics[width=0.55\linewidth]{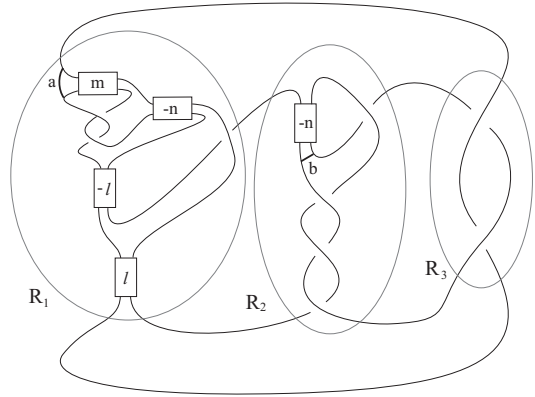}
\caption{$B(l, m, n, 0) + R(1)$ is a Montesinos link.}
\label{Blmn0+1Montesinos}
\end{center}
\end{figure}

\textsc{Proof of Lemma~\ref{Blmnp+R}.}
$(1)$ Figure~\ref{Blmn0+8unknot} shows that $B(l, m, n, 0) + R(\infty)$ is a trivial knot in $S^3$. 
In the Montesinos link of Figure~\ref{Blmn0+1Montesinos},
$R_1 = R(m, -2, -n, -l, -1, l, 0)$,
$R_2 = R(-n, -1, -3, 0)$, $R_3 = R(2, 0)$.
Assertion~(1)(ii) is obtained by computing
continued fractions.

\begin{figure}[h]
\begin{center}
\includegraphics[width=1.0\linewidth]{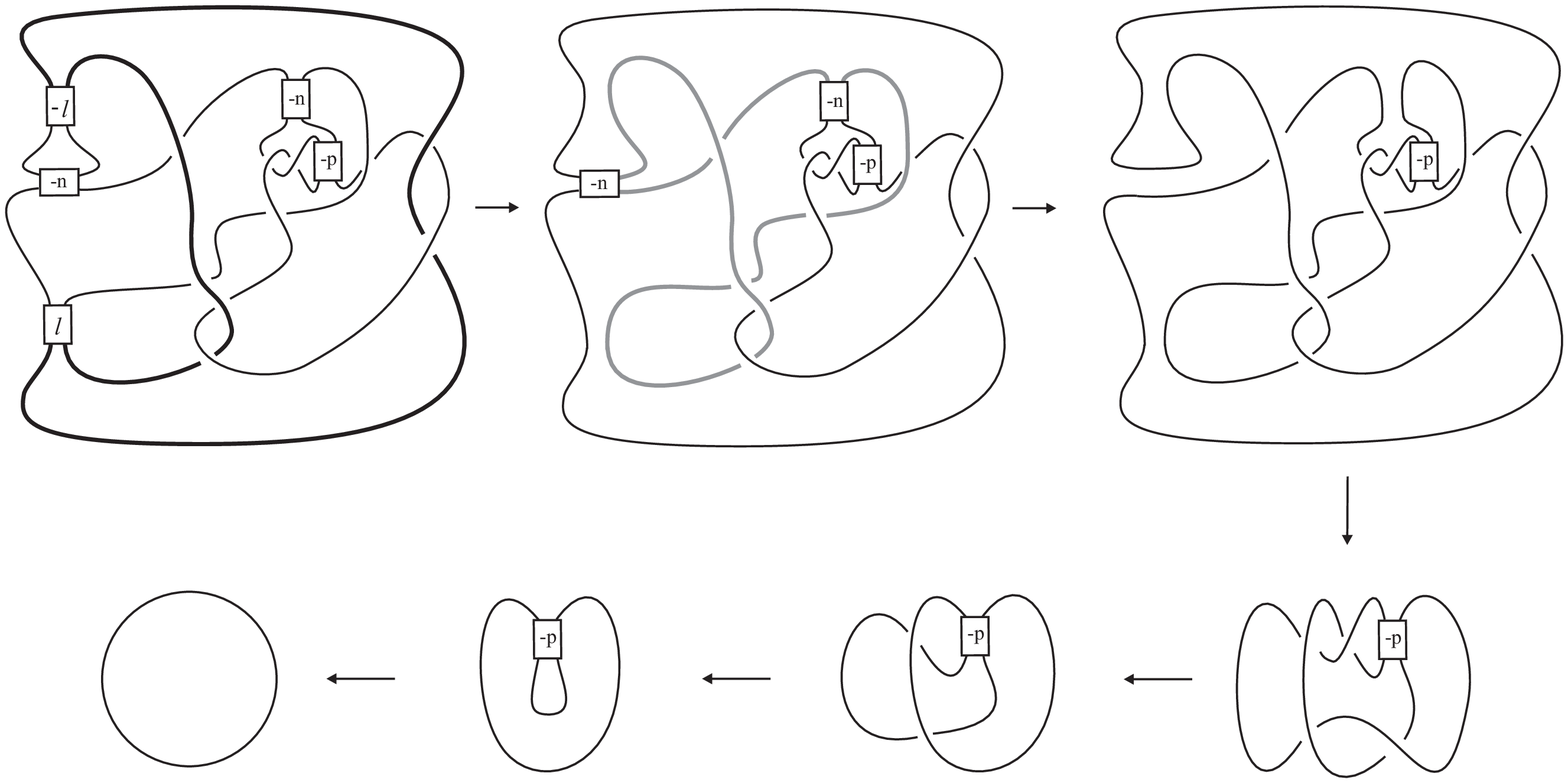}
\caption{ $B(l, 0, n, p) + R(\infty)$ is a trivial knot.}
\label{Bl0np+8unknot}
\end{center}
\end{figure}

\begin{figure}[h]
\begin{center}
\includegraphics[width=0.55\linewidth]{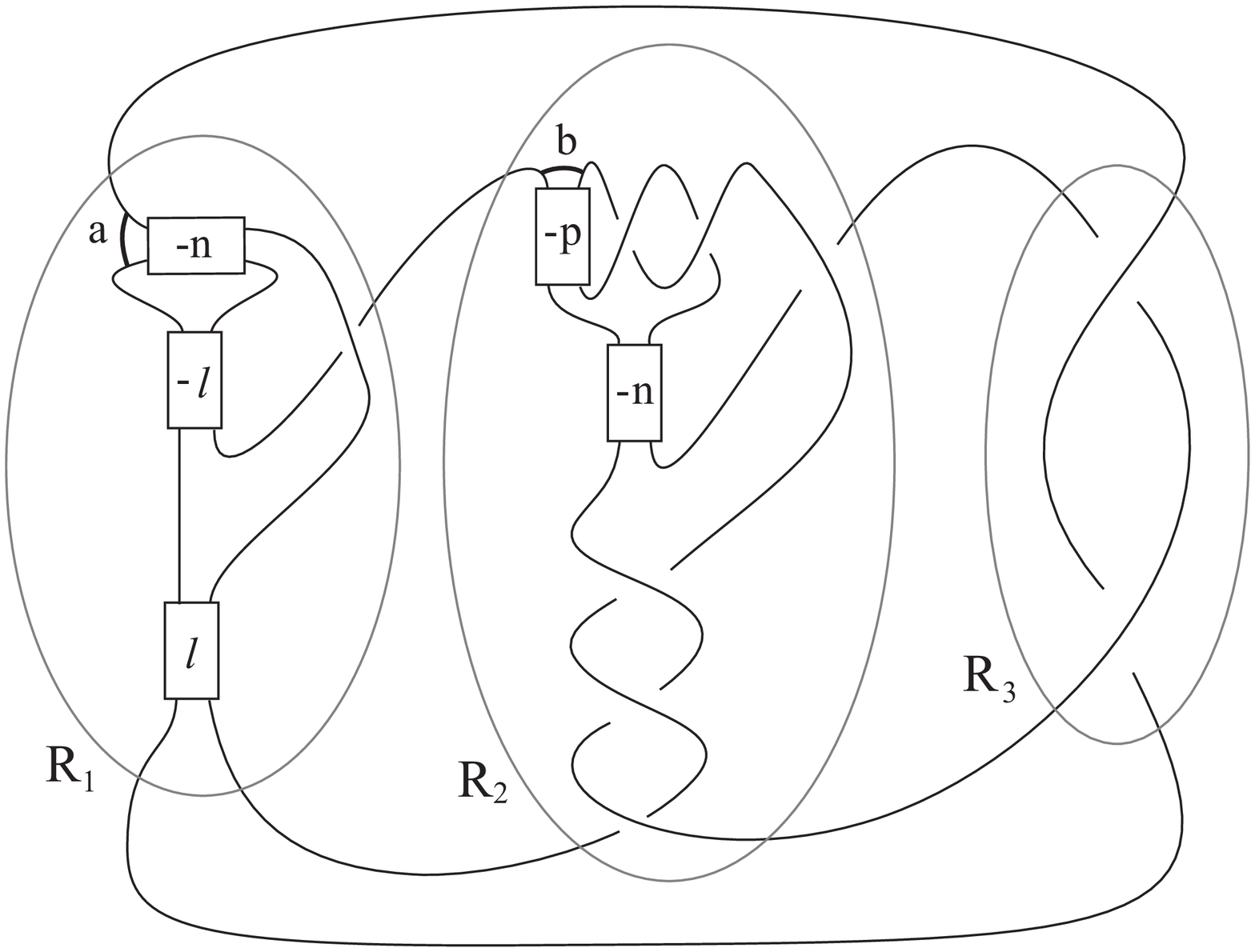}
\caption{$B(l, 0, n, p) + R(1)$ is a Montesinos link.}
\label{Bl0np+1Montesinos}
\end{center}
\end{figure}

$(2)$ 
Figure~\ref{Bl0np+8unknot} shows that $B(l, 0, n, p) + R(\infty)$ is a trivial knot in $S^3$. 
In the Montesinos link of Figure~\ref{Bl0np+1Montesinos},
$R_1 = R(-n, -l, -1, l, 0)$,
$R_2 = R(-p, 2, -n, -1, -3, 0)$, $R_3 = R(2, 0)$.
Assertion~(2)(ii) is obtained by computing
continued fractions.
\QED{Lemma~\ref{Blmnp+R}}

Let $K(l, m, n, p)$ be the covering knot 
of the trivializable tangle $B(l, m, n, p)$,
and $\gamma_{l,m,n,p}$ the covering slope corresponding to 
$1$--untangle surgery on $B(l, m, n, p) + R(\infty)$,
where $m$ or $p$ is $0$. 
As noticed in Section~\ref{section:covering}
the covering knot $K(l, m, n, p)$ is strongly invertible.
Then $0$--untangle surgery 
on $B(l, m, n, p) + R(\infty)$ 
corresponds to 
$(\gamma_{l, m, n, p} +1)$--surgery on $K(l, m, n, p)$
by Remark~\ref{n--framing}.
For brevity,
we often write $(K(l,m, n, p), \gamma)$ and
$(K(l, m, n, p), \gamma +1)$
for $(K(l,m,n,p), \gamma_{l,m,n,p})$
and $(K(l,m,n,p), \gamma_{l,m,n,p}+1)$,
respectively.
Lemma~\ref{Blmnp+R} shows that $(K(l,m, n, p), \gamma)$ is 
a Seifert fibered surgery,
and the resulting manifold is given in
Proposition~\ref{Klmnp} below.
Refer to the proof of Corollary~\ref{cor:location of Klmnp}
for the calculation of $\gamma_{l,m,n,p}$.

\begin{proposition}
\label{Klmnp}
\begin{enumerate}
\item
$K(l, m, n, 0)(\gamma_{l, m, n, 0})$ is a Seifert fiber space 

$\displaystyle
S^2(\frac{2lmn+lm-ln+2mn+3m-n-1}{2l^2mn+l^2m-l^2n+2lm-2m-l+1}, 
	-\frac{n+1}{4n+3}, \frac{1}{2})$. 
		
\smallskip
\item 
$K(l, 0, n, p)(\gamma_{l, 0, n, p})$ is a Seifert fiber space 

$\displaystyle
S^2(\frac{ln+n+1}{l^2n+l-1}, \frac{-2np+n-p+1}{8np-4n+2p-3},
\frac{1}{2})$. 
	
\end{enumerate}
Furthermore, 
$\gamma_{l, m, n, p} = 5 +l +n(l^2 +8l +12) +2n^2(l +2)^2
-m(2nl +4n +l +4)^2 -p(2nl +4n +2)^2$, where $m$ or $p$ is $0$.

\end{proposition}

\subsection{\bm{$(K(l,m,n,p), \gamma_{l,m,n,p} +1)$}
and primitive/Seifert positions}
\label{Non primitive/Seifert}

In this subsection,
we show that the Seifert fibered surgery $(K(l, m, n, p), \gamma_{l,m,n,p})$
($m = 0$ or $p = 0$) 
does not admit a primitive/Seifert position 
if $l,m,n,p$ satisfy more conditions (Proposition~\ref{nonPS}). 
For this purpose we study
$B(l, m, n, p) +R(0)$ and its $2$--fold branched cover.

\begin{figure}[h]
\begin{center}
\includegraphics[width=0.58\linewidth]{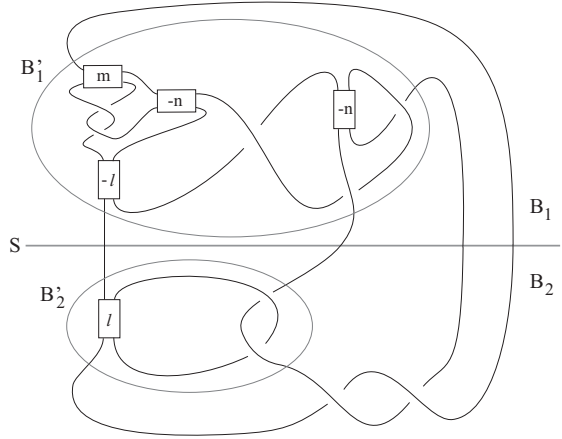}
\caption{$B(l, m, n, 0) + R(0)$ is a sum of Montesinos tangles.}
\label{Blmn0+0toroidal}
\end{center}
\end{figure}

\begin{figure}[h]
\begin{center}
\includegraphics[width=0.58\linewidth]{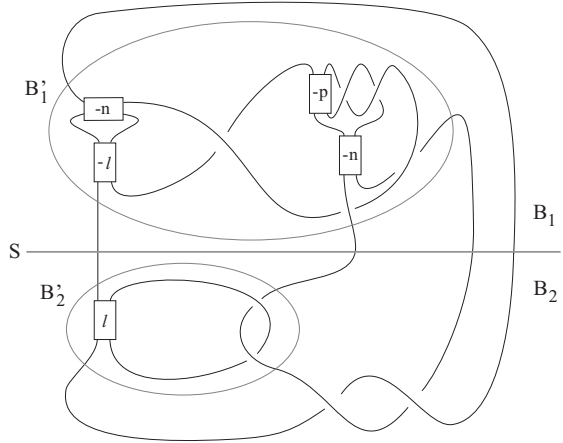}
\caption{$B(l, 0, n, p) + R(0)$ is a sum of Montesinos tangles.}
\label{Bl0np+0toroidal}
\end{center}
\end{figure}

Let $S$ be the $2$--sphere $\mbox{[plane]}\cup\{\infty\}$
intersecting $B(l, m, n, p) +R(0)$ in $4$ points
as in Figure~\ref{Blmn0+0toroidal} or \ref{Bl0np+0toroidal}
according as $p=0$ or $m=0$.
Let $B_1$ be the 3-ball bounded below by $S$,
and $B_2$ the 3-ball bounded above by $S$.

\begin{lemma}
\label{Blmnp+0}
\begin{enumerate}
\item
$B(l, m, n, 0) + R(0)$ is a sum of Montesinos tangles 

\smallskip\noindent
$\displaystyle
 M_T(-\frac{2lmn+lm-ln+2mn+3m-n-1}{2lmn+lm -ln+2m-1}, -\frac{n+1}{2n+1})$
 and $\displaystyle
M_T(\frac{1}{l}, -\frac{1}{2})$. 

\smallskip
\item
$B(l, 0, n, p) + R(0)$ is a sum of Montesinos tangles 

\smallskip\noindent
$\displaystyle
M_T(-\frac{ln+n+1}{ln+1}, -\frac{2np-n+p-1}{4np-2n-1})$
and 
$\displaystyle
M_T(\frac{1}{l}, -\frac{1}{2})$. 
\end{enumerate}
\end{lemma}

\noindent
\textsc{Proof of Lemma~\ref{Blmnp+0}.}
For brevity denote the knot or link $B(l, m, n, p) +R(0)$ by $K$,
where $m$ or $p$ is 0.

(1) Figure~\ref{Blmn0+0toroidal} shows that
the tangle $(B_1, B_1 \cap K)$ is pairwise homeomorphic
to $(B'_1, B'_1 \cap K) =
M_T(-\frac{2lmn+lm-ln+2mn+3m-n-1}{2lmn+lm -ln+2m-1}, -\frac{n+1}{2n+1})$  and
$(B_2, B_2 \cap K)$ is pairwise homeomorphic to
$(B'_2, B'_2\cap K) = M_T(\frac{1}{l}, -\frac{1}{2})$. 
It follows that $K$ is a sum of these two Montesinos tangles.

(2) Figure~\ref{Bl0np+0toroidal} shows that
$(B_1, B_1 \cap K)$ is pairwise homeomorphic to
$(B'_1, B'_1\cap K) =
M_T(-\frac{ln+n+1}{ln+1}, -\frac{2np-n+p-1}{4np-2n-1})$
and
$(B_2, B_2 \cap K)$ is pairwise homeomorphic to
$(B'_2, B'_2 \cap K) =
M_T(\frac{1}{l}, -\frac{1}{2})$.
\QED{Lemma~\ref{Blmnp+0}}

Let $\pi : K(l, m, n, p)(\gamma_{l,m,n,p}+1) \to S^3$ be the $2$--fold branched cover along $B(l, m, n, p) + R(0)$. 
Denote by $M_i$ the preimage $\pi^{-1}(B_i)$ $(i = 1, 2)$. 
By Lemma~\ref{Blmnp+0},
$M_1$ and $M_2$ are Seifert fiber spaces as described in
Proposition~\ref{KlmnpToroidal} below.
Proposition~\ref{KlmnpToroidal} shows that
the torus $\pi^{-1}(S)$
is an essential torus giving the torus decomposition
$K(l, m, n, p)(\gamma+1) =M_1 \cup M_2$.

\begin{proposition}
\label{KlmnpToroidal}
\begin{enumerate}
\item
Assume that 
$l \ne \pm 1, 0$, $n \ne 0, -1$, $(l, m, n) \ne  (-2, 0, 1)$, $(2, 1, -2)$. 
Then $\pi^{-1}(S)$ is a unique essential torus in
$K(l, m, n, 0)(\gamma_{l, m, n, 0}+1)$ up to isotopy,
and gives the torus decomposition with decomposing pieces 

\noindent
$\displaystyle M_1 = 
D^2(-\frac{2lmn+lm-ln+2mn+3m-n-1}{2lmn+lm -ln+2m-1}, -\frac{n+1}{2n+1})$ and

\noindent 
$\displaystyle M_2 = 
D^2(\frac{1}{l}, -\frac{1}{2})$. 

\smallskip
\item 
Assume that 
$l \ne \pm 1, 0$,\ $n \ne 0$, $(l, n) \ne (\pm 2, \mp 1)$,\ $(n, p) \ne (-1, 0)$,\ $(1, 1)$. 
Then $\pi^{-1}(S)$ is a unique essential torus in
$K(l, 0, n, p)(\gamma_{l, 0, n, p}+1)$ up to isotopy,
and gives the torus decomposition with decomposing pieces 

\noindent
$\displaystyle M_1 = D^2(-\frac{ln+n+1}{ln+1}, -\frac{2np-n+p-1}{4np-2n-1})$ and 
$\displaystyle M_2 =  D^2(\frac{1}{l}, -\frac{1}{2})$. 
\end{enumerate}
\end{proposition}

\noindent
\textsc{Proof of Proposition~\ref{KlmnpToroidal}.}
(1) We show that the Seifert fiber spaces $M_1$ and $M_2$ are boundary-irreducible
(i.e.\ $\pi^{-1}(S)$ is an essential torus),
and $M_1 \cup M_2$ is not a Seifert fiber space. 
Then the uniqueness of torus decomposition follows from 
\cite{JS, Jo}. 

The $2$--fold branched cover of a Montesinos tangle $M_T(\frac{p_1}{q_1}, \frac{p_2}{q_2})$ is a Seifert fiber space $D^2(\frac{p_1}{q_1},\ \frac{p_2}{q_2})$.
If $\frac{p_i}{q_i}$ is not an integer for $i=1, 2$,
then the Seifert fiber space is boundary-irreducible. 
We first show Claim~\ref{Seifert invariant1} below.

\begin{claim}
\label{Seifert invariant1}
$|2n +1| \ge 3$, $| l | \ge 2$, and $|2lmn+lm -ln+2m-1|\ge 2$.
\end{claim}

\noindent
\textsc{Proof of Claim~\ref{Seifert invariant1}.}
By the assumption of Proposition~\ref{KlmnpToroidal}(1),
$|2n+1| \ge 3$ and $|l| \ge 2$, 
so let us show $2lmn+lm -ln+2m-1 \ne 0, \pm 1$.  
Assume for a contradiction that
$2lmn +lm -ln +2m =\delta$ for some $\delta \in \{0, 1, 2\}$.
Then $m(2ln +l +2) =ln +\delta$. 
Since $|2n+1| \ge 3$, 
$2ln +l +2  =l(2n+1)+2 \ne 0$. 
If $m =0$, then $ln +\delta =0$.
This implies that $\delta \ne 0$ and
$(l, m, n) =(\pm1, 0, \mp\delta), (\pm\delta, 0, \mp1)$.
These are excluded by the assumption that
$l \ne 0, \pm1$, $n\ne -1$, $(l, m, n) \ne(-2, 0, 1)$.
Hence, we have $m \ne0$, so that
(*) $| 2ln + l +2 | \le | ln + \delta |$ holds.
Here we note that
$ln$, $ln +\delta$, $l(2n +1)$, and $l(2n+1) +2$
are of the same sign
because $\frac{2ln + l}{ln} = 2 + \frac{1}{n} > 0$, 
$| l(2n +1) | \ge 6$ and $| ln | \ge 2$ by the assumption. 
On the other hand, since $n \ne 0, -1$, 
it follows $|2n +1| \ge |n| +1$.
Hence, $| l(2n +1) | \ge | ln | + |l| \ge | ln | +2$,
so that $| l(2n +1) +2 | \ge | ln +\delta |$.
Therefore, (*) implies that the equality of (*) holds,
i.e.\ $l = 2, n = -2, \delta =0$.
It follows $(l, m, n) = (2, 1, -2)$,
which contradicts the assumption of Proposition~\ref{KlmnpToroidal}(1).
\QED{Claim~\ref{Seifert invariant1}}

It follows from Claim~\ref{Seifert invariant1} that
$M_2$ is boundary-irreducible.
Since $2lmn+lm-ln+2mn+3m-n-1$ and $2lmn+lm -ln+2m-1$,
and also $n+1$ and $2n+1$ are relatively prime,
Claim~\ref{Seifert invariant1} implies that none of
$\displaystyle \frac{2lmn+lm-ln+2mn+3m-n-1}{2lmn+lm -ln+2m-1}$ and 
$\displaystyle \frac{n+1}{2n+1}$ is an integer. 
Hence $M_1$ is also boundary-irreducible.
It follows that $\pi^{-1}(S)$ is an essential torus
in $K(l, m, 0, p)(\gamma_{l, m, 0, p}+1) = M_1\cup M_2$.

\begin{figure}[h]
\begin{center}
\includegraphics[width=1.0\linewidth]{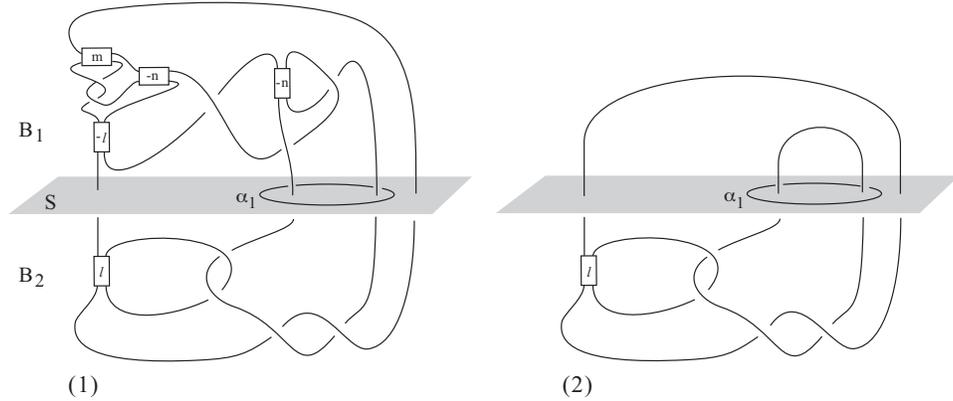}
\caption{$B(l,m,n,0)+R(0)$; $M_1 \in \mathcal{D}$, $M_2 \in \mathcal{L_K} \cup \mathcal{S_K}$.}
\label{Blmn0+0toroidal2_Lower}
\end{center}
\end{figure}

\begin{figure}[h]
\begin{center}
\includegraphics[width=1.0\linewidth]{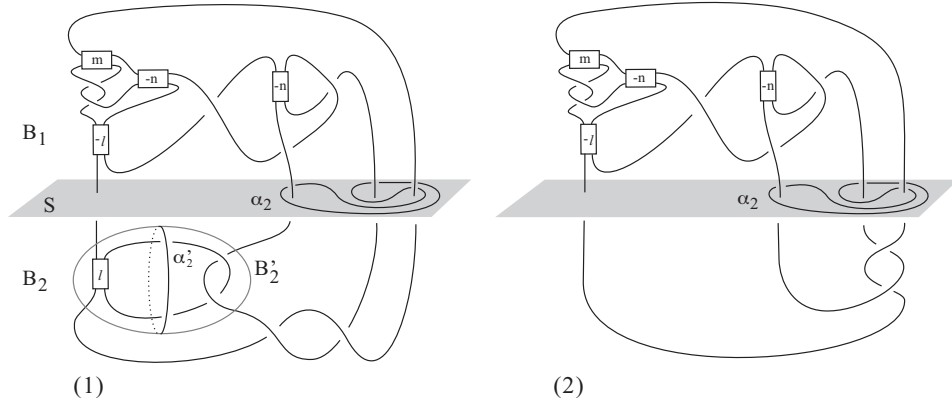}
\caption{$B(l,m,n,0)+R(0)$;
$M_1 \in \mathcal{L_K} \cup \mathcal{S_K}$, $M_2 \in \mathcal{D}$.}
\label{Blmn0+0toroidal3_Upper}
\end{center}
\end{figure}

\begin{figure}[h]
\begin{center}
\includegraphics[width=1.0\linewidth]{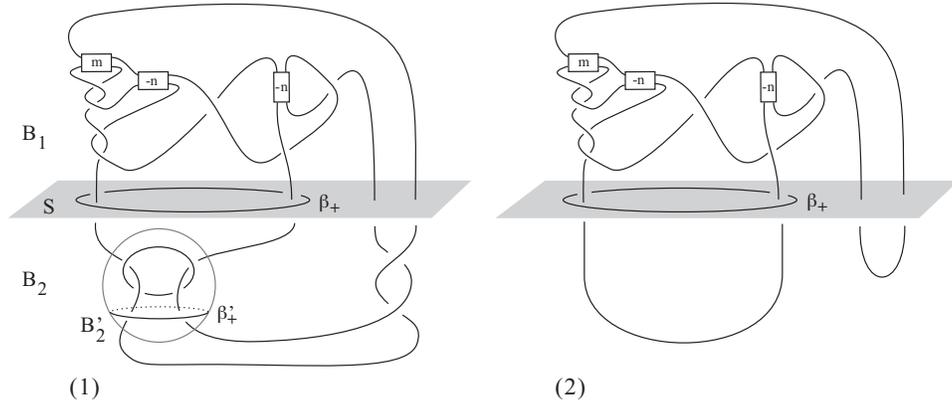}
\caption{$B(2,m,n,0)+R(0)$;
$M_1 \in \mathcal{S_K}$, $M_2 \in \mathcal{M}$.}
\label{B2mn0+0toroidal4_Upper}
\end{center}
\end{figure}

\begin{figure}[h]
\begin{center}
\includegraphics[width=1.0\linewidth]{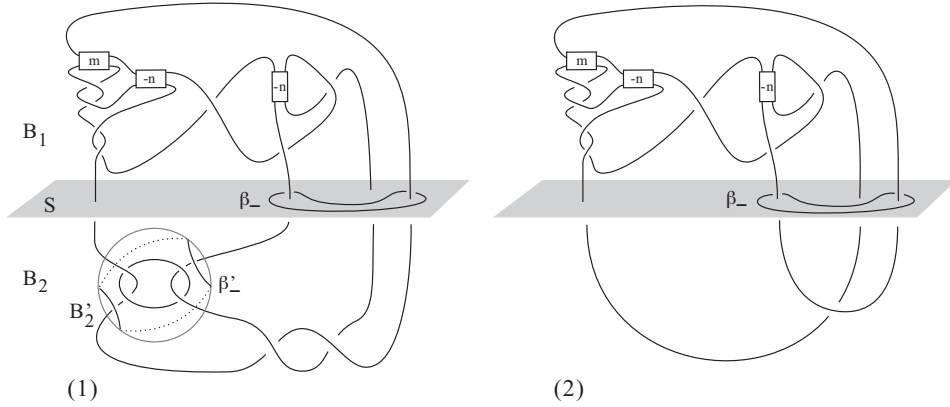}
\caption{$B(-2,m,n,0)+R(0)$;
$M_1 \in \mathcal{S_K}$, $M_2 \in \mathcal{M}$.}
\label{Bminus2mn0+0toroidal4_Upper}
\end{center}
\end{figure}

A Seifert fibration of $M_1$ is unique up to isotopy
(Lemma~\ref{lemma:fiber}),
so that its fiber in $\partial M_1 = \pi^{-1}(S)$
is isotopic to a lift of $\alpha_1 (\subset S)$
in Figure~\ref{Blmn0+0toroidal2_Lower}(1).
If $M_1 \cup M_2$ admits a Seifert fibration $\mathcal{F}$,
then the essential separating torus $\pi^{-1}(S)$ is a union
of fibers in $\mathcal{F}$.
This implies
$\mathcal{F}$ is an extension of fibrations of $M_1$ and $M_2$.
On the other hand, $M_2$ has a Seifert fibration
over the disk, and if $l = \pm 2$, then
$M_2$ has also a Seifert fibration over the M\"obius band.
In the former fibration,
lifts of $\alpha'_2 (\subset \partial B'_2)$ and thus
$\alpha_2 (\subset S)$ in Figure~\ref{Blmn0+0toroidal3_Upper}(1)
are fibers in $M_2$.
In the latter fibration with $l = 2$,
lifts of $\beta'_{+} (\subset \partial B'_2)$ and thus
$\beta_{+} (\subset S)$ in Figure~\ref{B2mn0+0toroidal4_Upper}(1)
are fibers in the Seifert fibration of $M_2$
over the M\"obius band (Lemma~\ref{lemma:fiber}).
If $l = -2$, lifts of $\beta'_{-} (\subset \partial B'_2)$ and thus
$\beta_{-} (\subset S)$ in Figure~\ref{Bminus2mn0+0toroidal4_Upper}(1)
are fibers in $M_2$. 
Since a lift of $\alpha_1$ is not isotopic in $\pi^{-1}(S)$
to any lift of $\alpha_2$ or $\beta_{\pm}$,
the fibrations of $M_1$ and $M_2$ do not match on their boundaries.
Hence, $M_1 \cup M_2$ is not a Seifert fiber space.

\begin{figure}[h]
\begin{center}
\includegraphics[width=1.0\linewidth]{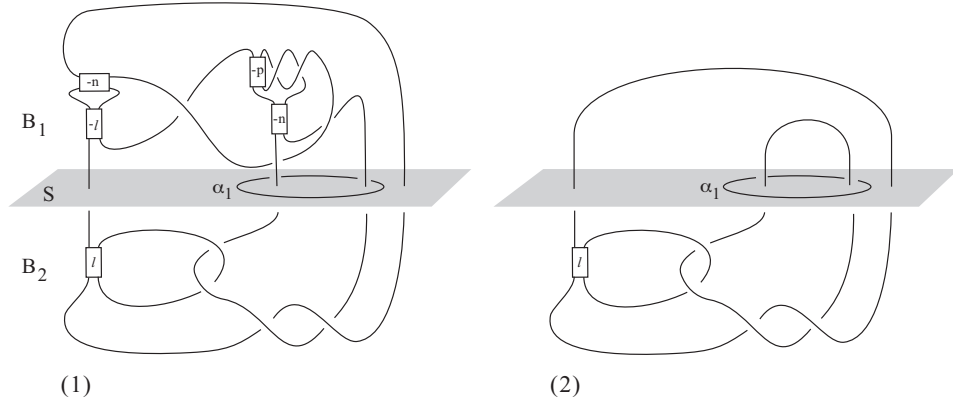}
\caption{$B(l,0,n,p)+R(0)$;
$M_1 \in \mathcal{D}$, $M_2 \in \mathcal{L_K} \cup \mathcal{S_K}$.}
\label{Bl0np+0toroidal2_Lower}
\end{center}
\end{figure}

\begin{figure}[h]
\begin{center}
\includegraphics[width=1.0\linewidth]{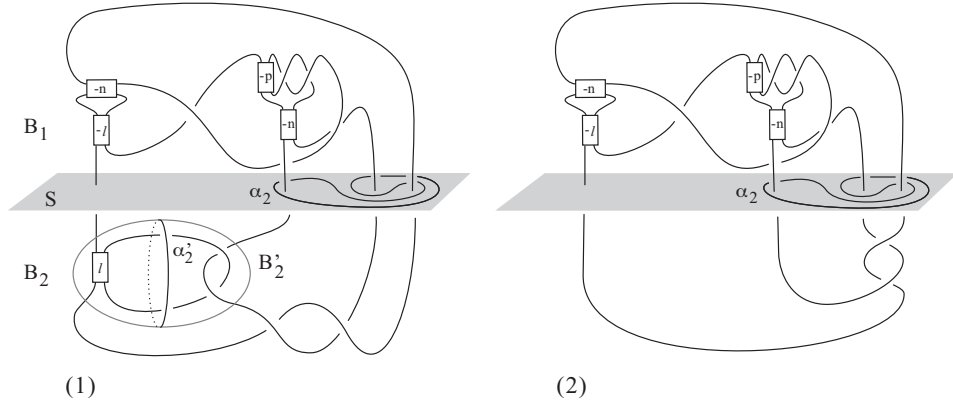}
\caption{$B(l,0,n,p)+R(0)$;
$M_1 \in \mathcal{L_K} \cup \mathcal{S_K}$, $M_2 \in \mathcal{D}$}
\label{Bl0np+0toroidal3_Upper}
\end{center}
\end{figure}

\begin{figure}[h]
\begin{center}
\includegraphics[width=1.0\linewidth]{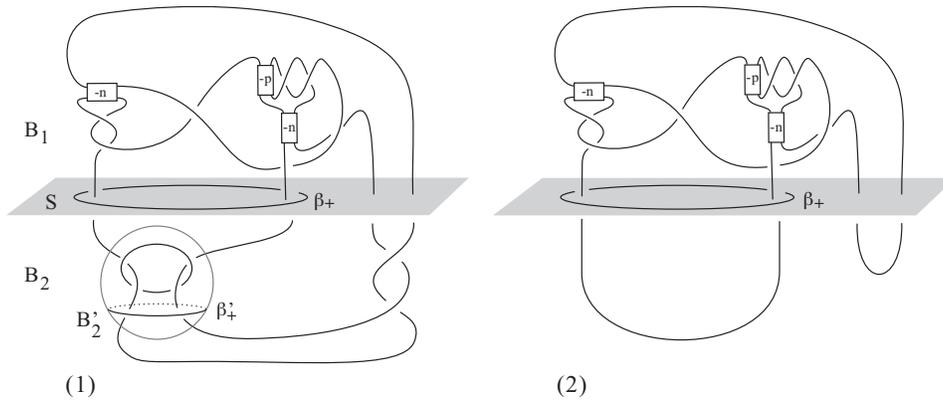}
\caption{$B(2,0,n,p)+R(0)$;
$M_1 \in \mathcal{S_K}$, $M_2 \in \mathcal{M}$.}
\label{B20np+0toroidal4_Upper}
\end{center}
\end{figure}

\begin{figure}[h]
\begin{center}
\includegraphics[width=1.0\linewidth]{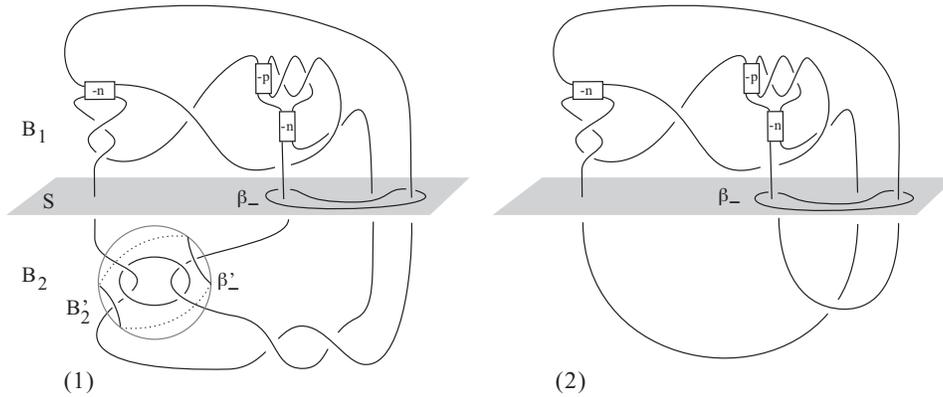}
\caption{$B(-2,0,n,p)+R(0)$;
$M_1 \in \mathcal{S_K}$, $M_2 \in \mathcal{M}$.}
\label{Bminus20np+0toroidal4_Upper}
\end{center}
\end{figure}

(2) The assumption on $l, n, p$ in
Proposition~\ref{KlmnpToroidal}(2)  
assures that $| l |\ge2$, $|ln+1| \ge 2$, and $|4np-2n-1|\ge 3$,
hence $M_1$ and $M_2$ 
are boundary-irreducible Seifert fiber spaces.
As in (1) we see that 
Seifert fibrations of $M_1$ and $M_2$
do not match on their boundaries;
a lift of $\alpha_1 (\subset S)$ is a fiber in $M_1$
in Figures~\ref{Bl0np+0toroidal2_Lower}(1), and
a lift of $\alpha_2$, $\beta_{+}$, or $\beta_{-}$$(\subset S)$
in Figures~\ref{Bl0np+0toroidal3_Upper}(1),
\ref{B20np+0toroidal4_Upper}(1), \ref{Bminus20np+0toroidal4_Upper}(1)
 is a fiber in $M_2$.
It follows that $M_1 \cup M_2$ is not a Seifert fiber space,
as desired.
\QED{Proposition~\ref{KlmnpToroidal}}

Now we are ready to state and prove the main result of
this section:
under some conditions on $l, m, n, p$ slightly stronger
than those in Proposition~\ref{KlmnpToroidal},
the Seifert fibered surgery $(K(l, m, n, p), \gamma)$
does not have a primitive/Seifert position.

\begin{proposition}
\label{nonPS}
\begin{enumerate}
\item 
Assume that $l \ne \pm 1, 0$, $n \ne 0, -1$, 
$(l, m) \ne (-2, 0)$, $(-2, 2)$,  
$(l, m, n) \ne  (2, 1, -2)$.  
Then the Seifert fibered surgery $(K(l, m, n, 0), \gamma_{l,m,n,0})$ does not have a primitive/Seifert position. 

\item
Assume that $l \ne \pm 1, 0$, $n \ne  0$, $(l, n) \ne (\pm 2, \mp 1)$,
$(l, p) \ne (-2, 0)$, $(-2, 2)$, 
$(n, p) \ne (-1, 0)$, $(1,1)$. 
Then the Seifert fibered surgery $(K(l, 0, n, p), \gamma_{l, 0, n, p})$ does not have a primitive/Seifert position.  
\end{enumerate}
\end{proposition}

\noindent
\textsc{Proof of Proposition~\ref{nonPS}.} 
If a Seifert fibered surgery on a knot $K$ has a primitive/Seifert position,  
then $K$ has a tunnel number one \cite[2.3]{D}. 
Thus the result follows from Proposition~\ref{tunnel2} below. 
\QED{Proposition~\ref{nonPS}}

\begin{proposition}
\label{tunnel2}
\begin{enumerate}
\item 
Assume that $l, m, n$ satisfy the condition in
Proposition~\ref{nonPS}(1). 
Then the tunnel number of $K(l, m, n, 0)$ is two.  

\item
Assume that $l, n, p$ satisfy the condition in
Proposition~\ref{nonPS}(2). 
Then the tunnel number of $K(l, 0, n, p)$ is two.   
\end{enumerate}
\end{proposition}

\noindent
\textsc{Proof of Proposition~\ref{tunnel2}.}
Assume for a contradiction that $K(l, m, n, p)$ ($m$ or $p$ is $0$) has tunnel number one. 
Then $K(l, m, n, p)(r)$ admits a genus two Heegaard splitting
for any slope $r$.
Let $r =\gamma_{l, m, n, p} +1$.
We see from Proposition~\ref{KlmnpToroidal} that
$K(l, m, n, p)(r)$ contains
a unique essential torus up to isotopy,
which decompose $K(l, m, n, p)(r)$ into
two Seifert fiber spaces $M_1,M_2$.

Kobayashi \cite[Theorem]{Koba} classifies toroidal $3$--manifolds
with genus two Heegaard splittings.
In our setting,
case~(i), (ii), or (iii) in Theorem in \cite{Koba} holds,
and we see that
 $K(l, m, n, p)(r)$ is obtained by gluing
pieces in classes $\mathcal{D}, \mathcal{M}, \mathcal{S_{K}}, 
\mathcal{L_{K}}$ defined below.
The class $\mathcal{D}$ is the set of Seifert fiber spaces
over the disk with two exceptional fibers, and
$\mathcal{M}$ is the set of Seifert fiber spaces over the M\"obius band with no exceptional fiber. 
The class $\mathcal{S_K}$ is the set of the exteriors of
2-bridge knots in the $3$--sphere, and
$\mathcal{L_K}$ is the set of the exteriors of $1$--bridge knots in lens spaces;
a $1$--bridge knot $K$ in a lens space $M$ is 
a knot intersecting a genus one Heegaard surface $F$ of
$M =V_1 \cup_{F} V_2$ in two points
such that $K \cap V_i$ is a boundary parallel arc in $V_i$
for $i=1, 2$.
Applying Theorem in \cite{Koba},
we obtain the following lemma on $M_1, M_2$.

\begin{lemma}[\cite{Koba}]
\label{Heegaard}
There are the following four possibilities on $M_1$ and $M_2$.
\begin{enumerate}
\item
$M_1 \in \mathcal{D}$, $M_2 \in \mathcal{L_K} \cup \mathcal{S_K}$

\item
$M_1 \in \mathcal{L_K} \cup \mathcal{S_K}$, $M_2 \in \mathcal{D}$ 

\item
$M_1 \in \mathcal{M}$, $M_2 \in \mathcal{S}_K$

\item
$M_1 \in  \mathcal{S_K}$, $M_2 \in \mathcal{M}$
\end{enumerate}
Furthermore,
in $(1)$, $(3)$ a regular fiber of $M_1$ is identified with a meridian of $M_2$, and in $(2)$, $(4)$
a regular fiber of $M_2$ is identified with a meridian of $M_1$.
\end{lemma}

We first assume $p =0$.
Let us derive a contradiction
in each case of Lemma~\ref{Heegaard}.
For brevity denote $K = B(l, m, n, 0) +R(0)$.

Assume that case~(1) in Lemma~\ref{Heegaard} occurs.
A lift of the simple closed curve $\alpha_1 (\subset S)$
in Figure~\ref{Blmn0+0toroidal2_Lower}(1)
is a regular fiber of $M_1$ contained in the torus
$\partial M_1 = \pi^{-1}(S)$.
Replace the tangle $(B_1, B_1 \cap K)$
in Figure~\ref{Blmn0+0toroidal2_Lower}(1) with
a trivial tangle $(B, t)$ in which
 $\alpha_1$ bounds a disk in $B -t$.
Then the $2$--fold branched cover of $(B, t) \cup (B_2, B_2 \cap K)$
is obtained by gluing a solid torus $V$ to
$M_2$ along its boundary $\partial M_2$ so that
a meridian of $V$ is identified with
a lift of $\alpha_1 (\subset S)$.
Since the lift is a meridian of
the exterior $M_2$ of a knot in a lens space or the $3$--sphere
by Lemma~\ref{Heegaard},
$V \cup M_2$ is either a lens space or the $3$--sphere.
However, $(B, t) \cup (B_2, B_2 \cap K)$ is,
as depicted in Figure~\ref{Blmn0+0toroidal2_Lower}(2),
the Montesinos link $M(\frac{1}{l}, -\frac{1}{2}, \frac{1}{2})$,
where $| l | \ge 2$.
It follows that $V \cup M_2$ is a Seifert fiber space 
with three exceptional fibers,
which is not a lens space or the $3$--sphere. 
Hence, (1) in Lemma~\ref{Heegaard} does not occur.

Assume that case~(2) in Lemma~\ref{Heegaard} occurs.
Then, a lift of the simple closed curve $\alpha_2 (\subset S)$
in Figure~\ref{Blmn0+0toroidal3_Upper}(1) is a regular fiber
of the Seifert fibration of $M_2$ over the disk.
Replacing $(B_2, B_2 \cap K)$ in
Figure~\ref{Blmn0+0toroidal3_Upper}(1) with a trivial tangle $(B, t)$
in which $\alpha_2$ bounds a disk in $B -t$,
we obtain the Montesinos link
$$M(-\frac{2lmn+lm-ln+2mn+3m-n-1}{2lmn+lm -ln+2m-1},
-\frac{n+1}{2n+1}, \frac{1}{2})$$
as depicted in
 Figure~\ref{Blmn0+0toroidal3_Upper}(2);
note $| 2lmn+lm -ln+2m-1 |\ge 2$, $|2n+1| \ge 3$
by Claim~\ref{Seifert invariant1}.
The $2$--fold branched cover along this Montesinos link
is a Seifert fiber space with three exceptional fibers.
However, since $M_1 \in \mathcal{L_K}\cup\mathcal{S_K}$
and a lift of $\alpha_2$ is a meridian of the knot exterior $M_1$,
by the same arguments as in case~(1) the $2$--fold branched
cover is a lens space or the $3$--sphere.
This is a contradiction.

Assume that case~(3) in Lemma~\ref{Heegaard} occurs.
By Proposition~\ref{KlmnpToroidal}(1)
$M_1$ has a Seifert fibration over the disk with
an exceptional fiber of index $|2n +1|$, an odd integer.
Hence, $M_1 \not\in \mathcal{M}$ by Lemma~\ref{lemma:fiber},
a contradiction.

Assume that case~(4) in Lemma~\ref{Heegaard} occurs.
Then $l = \pm 2$, 
and $M_2 = D^2(\frac{1}{l}, -\frac{1}{2})$ has a Seifert fibration over the M\"obius band. 
In this Seifert fibration, 
a regular fiber in $\partial M_2$ is isotopic to
a lift of $\beta_{+} (\subset S)$
given in Figure~\ref{B2mn0+0toroidal4_Upper}(1) if $l =2$, and
a lift of $\beta_{-} (\subset S)$
in Figure~\ref{Bminus2mn0+0toroidal4_Upper}(1) if $l = -2$.

Replacing $(B_2, B_2 \cap K)$ with a trivial tangle $(B, t)$
in Figure~\ref{B2mn0+0toroidal4_Upper}(1)
(resp.\ \ref{Bminus2mn0+0toroidal4_Upper}(1)) 
with the rational tangle $(B, t)$
in which $\beta_{+}$ (resp.\ $\beta_{-}$) bounds a disk in $B -t$,
we obtain the Montesinos link
depicted in Figure~\ref{B2mn0+0toroidal4_Upper}(2)
(resp.\ \ref{Bminus2mn0+0toroidal4_Upper}(2)).
The $2$--fold cover branched along this link is
the Seifert fiber space $S_{\pm l, m, n}$ as follows:

\medskip
\noindent
$\displaystyle
S_{2,m,n} = S^2(\frac{-6mn-5m+3n+1}{4mn+4m-2n-1}, \frac{-n-1}{2n+1})$
if $l = 2$, and

\medskip\noindent
$\displaystyle
S_{-2, m,n} = S^2(\frac{-2mn+m+n-1}{4mn-2n+1}, \frac{n}{2n+1})$
if $l = -2$.

\medskip\noindent
Since $M_1 \in \mathcal{S_K}$
and a lift of $\beta_{\pm}$ is a meridian of the exterior $M_1$ 
of a knot in $S^3$,
by the same argument as in case~(1) the $2$--fold branched
cover $S_{\pm2, m, n}$ is $S^3$.

\begin{claim}
\label{homology1}
\begin{enumerate}
\item
If $S_{2, m, n}$ is $S^3$, then $n = -1$.
\item
$S_{-2, m, n}$ is $S^3$ if and only if $m = 0, 2$. 
\end{enumerate}
\end{claim}

\noindent
\textsc{Proof of Claim~\ref{homology1}.}
(1) 
Note that the first homology group of $S_{2,m,n}$ is 
the cyclic group of order
$|(-6mn-5m+3n+1)(2n+1) + (-n-1)(4mn+4m-2n-1)| = 
| 16mn^2 + 24mn -8n^2 + 9m -8n -2 |$.  
Assume that the order equals $1$.
Then, $(16n^2+24n+9)m = 8n^2+8n+ \delta$,
where $\delta$ is $1$ or $3$.
Note that $16n^2+24n+9 > 0$ and $8n^2 +8n +\delta > 0$ 
for any integer $n$.
It follows that $16n^2+24n+9 \le 8n^2+8n+\delta$ for some integer $n$.
This inequality has the only integral solution $n =-1$.

(2) 
The first homology of $S_{-2,m,n}$ has order 
$|(-2mn+m+n-1)(2n+1) + n(4mn-2n+1)| = |m-1|$, 
which is $1$ if and only if $m = 0, 2$. 
\QED{Claim~\ref{homology1}}

Claim~\ref{homology1}, together with the assumption in Proposition~\ref{nonPS}(1), 
shows that $S_{\pm 2, m, n}$ is not $S^3$. 
Thus case~(4) in Lemma~\ref{Heegaard} does not occur.
Hence, the tunnel number of $K(l, m, n, 0)$ is greater than one.

Assume $m =0$.
Using the above arguments for $p=0$,
we prove that cases~(1), (2), (3), (4) in Lemma~\ref{Heegaard}
do not occur.
Follow the proof for $p=0$ with
Figures~\ref{Blmn0+0toroidal2_Lower}, \ref{Blmn0+0toroidal3_Upper}, 
\ref{B2mn0+0toroidal4_Upper} and \ref{Bminus2mn0+0toroidal4_Upper}
replaced by 
Figures~\ref{Bl0np+0toroidal2_Lower}, \ref{Bl0np+0toroidal3_Upper}, 
\ref{B20np+0toroidal4_Upper} and \ref{Bminus20np+0toroidal4_Upper}, respectively. 
Cases~(1) and (2) in Lemma~\ref{Heegaard} do not occur, 
because the Montesinos links in
Figures~\ref{Bl0np+0toroidal2_Lower}(2),
\ref{Bl0np+0toroidal3_Upper}(2)
consist of three rational tangles and their $2$--fold branched cover
cannot be a lens space or the 3-sphere.

Case~(3) in Lemma~\ref{Heegaard} does not occur.
This is because 
the Seifert fibration of $M_1$ over the disk 
(Proposition~\ref{KlmnpToroidal}(2)) contains
an exceptional fiber of index $|4np -2n-1|$, an odd integer,
and thus $M_1 \not \in \mathcal{M}$.

Assume case~(4) in Lemma~\ref{Heegaard} occurs;
then $l=\pm2$ and
the $2$--fold cover branched along the link in
Figure~\ref{B20np+0toroidal4_Upper}(2) or
\ref{Bminus20np+0toroidal4_Upper}(2)
is the $3$--sphere.
If $l = 2$,
then the Montesinos link in Figure~\ref{B20np+0toroidal4_Upper}(2)
is $\displaystyle
M(-\frac{3n +1}{2n +1}, -\frac{2np -n +p-1}{4np -2n -1})$,
and the first homology group of
the $2$--fold branched cover along this link has
order $|16n^2p-8n^2+8np-8n +p-2|$. 
If the order is $1$, then
$(16n^2 +8n +1)p = 8n^2+8n +\delta$ where $\delta$ is $1$ or $3$.
Since $16n^2 +8n +1 > 0$ and $8n^2+8n +\delta >0$
for any integer $n$,
we obtain $16n^2 +8n +1 \le 8n^2+8n +\delta$.
Then $n =0$.
This contradicts the assumption $n \ne 0$. 
If $l =-2$,
then the Montesinos link in
Figure~\ref{Bminus20np+0toroidal4_Upper}(2)
is $\displaystyle M(\frac{n -1}{-2n +1}, \frac{2np -n -p}{4np -2n -1})$,
and the first homology group of
the $2$--fold branched cover along this link has
order $|p-1|$. 
By the arguments similar to above,
we see that if the order is $1$, 
then $p = 0$ or $2$.  
This contradicts the assumption $(l, p) \ne (-2, 0),\ (-2, 2)$. 
So, case~(4) does not occur and thus
$K(l, 0, m, p)$ has tunnel number greater than one.

\begin{figure}[h]
\begin{center}
\includegraphics[width=0.6\linewidth]{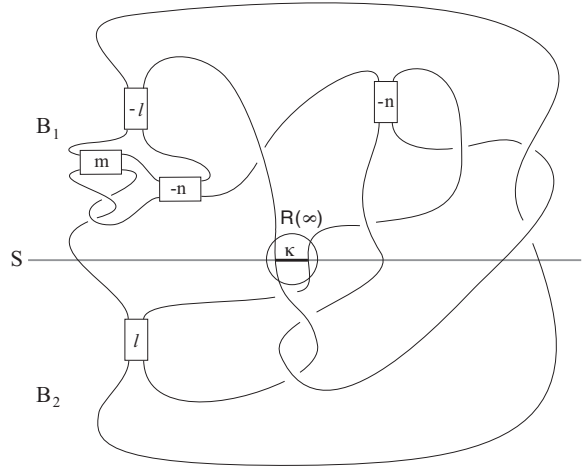}
\caption{Both sides of $S$ are $3$--string trivial tangles.}
\label{tunnel2Blmnp+8}
\end{center}
\end{figure}

Let us show that the tunnel number of $K(l, m, n, p)$
equals 2, where $m$ or $p$ is 0.
For brevity denote $K = B(l, m, n, p) + R(\infty)$.
Let $S$ be the $2$--sphere $\mbox{[plane]}\cup \{\infty\}$
given in Figure~\ref{tunnel2Blmnp+8},
and $B_1$ (resp.\ $B_2$) the $3$--balls bounded
below (resp.\ above) by $S$.
Then, each $(B_i, B_i \cap K)$ is a $3$--string tangle pairwise
homeomorphic to $(D^2 \times I, \{ x_1, x_2, x_3 \}\times I)$,
where $x_i$ are distinct points in $D^2$.
Let $\pi_{\infty} : S^3 \to S^3$ be the $2$--fold cover branched along
the trivial knot $K$.
The preimage $\pi_{\infty}^{-1}(B_i)$ is a genus $2$ handlebody
for $i=1, 2$,
and thus $\pi_{\infty}^{-1}(S)$ is
a genus 2 Heegaard surface of $S^3$.
Since $S$ contains the spanning arc $\kappa$ for $R(\infty)$
as described in Figure~\ref{tunnel2Blmnp+8},
the covering knot $K(l, m, n, p) =\pi_{\infty}^{-1}(\kappa)$
is contained in the genus 2 Heegaard surface $\pi_{\infty}^{-1}(S)$.
Then, by \cite[Fact on p.138]{Mor}
the tunnel number of $K(l, m, n, p)$ is at most $2$,
and so equals $2$. 
\QED{Proposition~\ref{tunnel2}}

\subsection{Hyperbolicity of \bm{$K(l,m,n,p)$}}
\label{subsection:hyperbolicity}

We use the following lemma to detect hyperbolicity of knots.

\begin{lemma}
\label{lemma:hyperbolic}
Suppose that an $r$--surgery on a knot $K$ in $S^3$ yields
a $3$--manifold $K(r)$ containing a separating incompressible torus.
Suppose further that no separating incompressible torus
in $K(r)$ is disjoint from the dual knot $K^*$ of $K$,
i.e.\ the core of the filled solid torus in $K(r)$.
Then, $K$ is a hyperbolic knot.
\end{lemma}

\noindent
\textsc{Proof of Lemma~\ref{lemma:hyperbolic}.}
If $K$ is not a hyperbolic knot, then
it is either a torus knot or a satellite knot.
If the result of a surgery on a torus knot
contains an incompressible torus, 
then the surgery is longitudinal \cite[VI.\ Example]{J}
and the torus is non-separating. 
Hence, $K$ is not a torus knot
because $K(r)$ contains a separating incompressible torus
by the assumption.
Now assume that $K$ is a satellite knot with
a companion knot $k$.
Then $K$ has a companion knot $k$ which is either
a torus knot or a hyperbolic knot; 
$K$ is contained in a tubular neighborhood $V$ of $k$.
Since the separating torus $\partial V$ is disjoint from $K^*$ in $K(r)$,
the assumption of the lemma implies that
$\partial V$ compresses after the $r$--surgery along $K(\subset V)$.
By \cite{Ga1}
$K$ is a $0$ or $1$--bridge braid
in $V$ and winds $w (\ge 2)$ times in $V$. 
It follows that $K(r) = k(\frac{m}{nw^2})$ \cite{Go1},
where $r = \frac{m}{n}$, and $m$ and $w^2$ are relatively prime. 
Since $k(\frac{m}{n w^2})$ contains a separating incompressible torus
by the assumption, 
$k$ is not a torus knot and thus a hyperbolic knot. 
However, \cite{GL4} shows that
if $k(\frac{m}{n w^2})$ is toroidal,
then $|n w^2| \le 2$, a contradiction. 
Hence, $K$ is not a satellite but a hyperbolic knot.
\QED{Lemma~\ref{lemma:hyperbolic}}

\begin{proposition}
\label{hyperbolic}
The covering knot $K(l, m, n, p)$ $($$m$ or $p$ is $0$$)$
is a hyperbolic knot
if $l, m, n, p$ satisfy the condition in
Proposition~\ref{KlmnpToroidal}. 
\end{proposition}

\noindent
\textsc{Proof of Proposition~\ref{hyperbolic}.}
We prove that $K =K(l, m, n, p)$ and $r = \gamma_{l,m,n,p} +1$
satisfy the assumption of Lemma~\ref{lemma:hyperbolic}.
Proposition~\ref{KlmnpToroidal} shows that
$\pi^{-1}(S)$ is a unique incompressible torus
in $K(l, m, n, p)(\gamma_{l,m,n,p}+1)$ up to isotopy,
where $\pi: K(l,m,n,p)(\gamma_{l,m,n,p} +1) \to S^3$ is
the $2$--fold cover branched along $B(l,m,n,p)+R(0)$. 
Thus, Claim~\ref{intersection} below shows that
the assumption is satisfied when $p =0$,
so that $K(l, m, n, 0)$ is a hyperbolic knot.
Claim~\ref{intersection}
with $l, m, n, 0$ replaced by $l, 0, n, p$ also holds,
and implies that $K(l, 0, n, p)$ is hyperbolic.
\QED{Proposition~\ref{hyperbolic}}

\begin{claim}
\label{intersection}
In $K(l, m, n, 0)(\gamma_{l,m,n,0} +1)$,
the incompressible torus $\pi^{-1}(S)$ intersects
the dual knot $K^*$ of $K(l, m, n, 0)$ minimally in two points. 
\end{claim}

\begin{figure}[h]
\begin{center}
\includegraphics[width=1.0\linewidth]{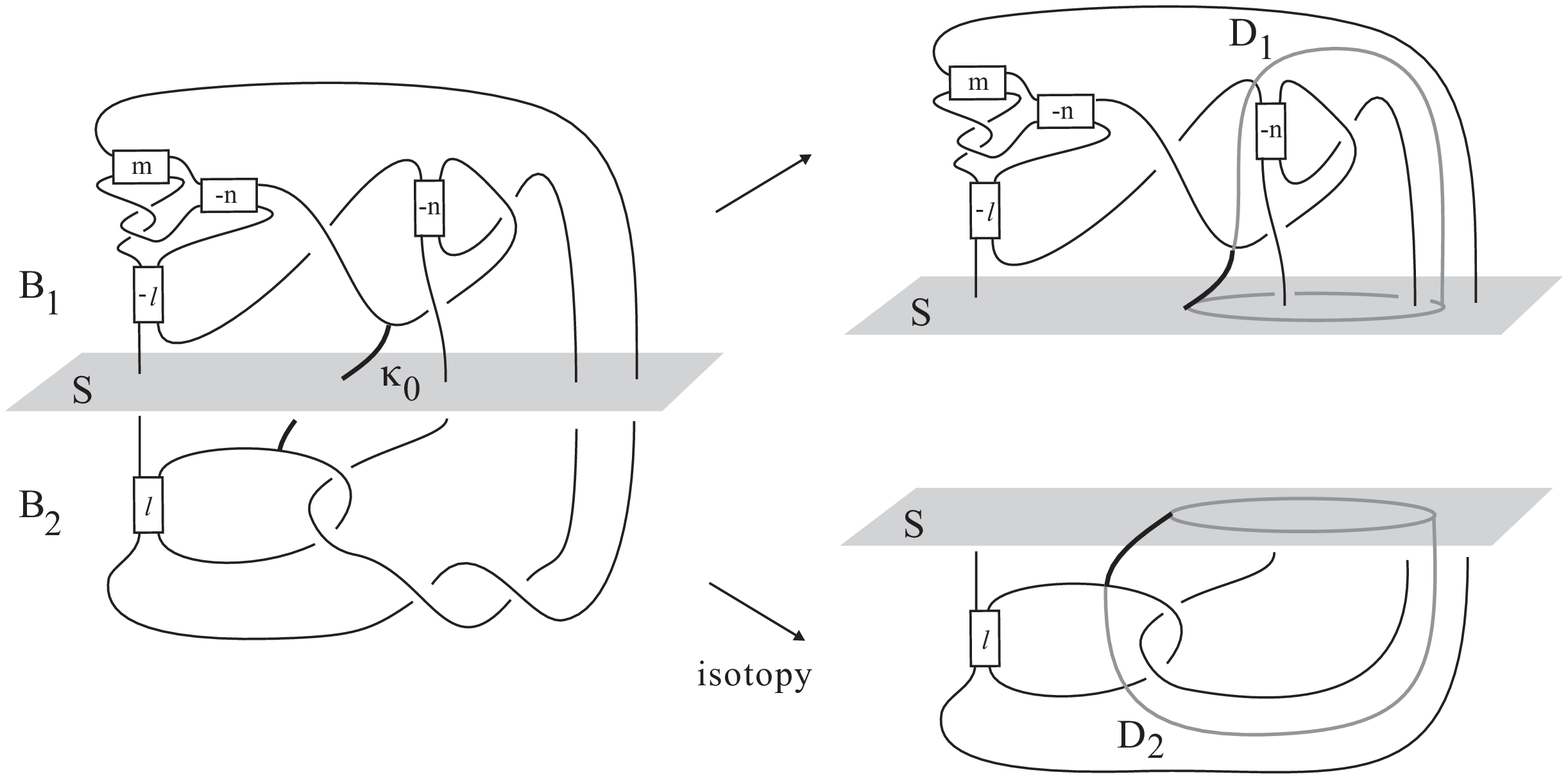}
\caption{$B(l, m, n, 0)+R(0)$.}
\label{hyperbolicityKlmnp}
\end{center}
\end{figure}

\noindent
\textsc{Proof of Claim~\ref{intersection}.}
The arc $\kappa_0$ in Figure~\ref{hyperbolicityKlmnp}
is a spanning arc of $R(0)$.
The preimage $\pi^{-1}(\kappa_0)$ is the dual knot $K^*$.
In Figure~\ref{hyperbolicityKlmnp},
$D_i$ is a disk which 
contains $\kappa_0 \cap B_i$ and
splits the tangle $(B_i, B_i \cap B(l,m,n,0)+R(0))$
into two nontrivial tangles.
It follows that $A_i = \pi^{-1}(D_i)$ is an essential annulus
in the Seifert fiber space $M_i = \pi^{-1}(B_i)$ over the disk with 
two exceptional fibers. 
Furthermore,
the arc $\pi^{-1}(\kappa_0 \cap B_i) = K^* \cap M_i$
is an essential arc in the annulus $A_i$. 
Then, the desired result follows from
the argument in \cite[Example 1.4]{EM0}.
\QED{Claim~\ref{intersection}}

\section{Locating {$(K(l,m,n,p), \gamma_{l,m,n,p})$}
in the Seifert Surgery Network}
\label{Seifert Surgery Network}

In Subsection~\ref{subsection:seiferters and tangles},
we review a method of finding seiferters for Seifert surgeries
obtained by untangle surgeries.
In Subsection~\ref{Location},
we find  seiferters for the Seifert fibered surgeries
$(K(l,m,n,p), \gamma_{l,m,n,p})$,
and an explicit path to a Seifert surgery on a trefoil knot.

\subsection{Seiferters and tangles}
\label{subsection:seiferters and tangles}

Assume that
a tangle $(B, t)$, where $B$ is the complement of
the unit $3$--ball in $S^3$,
satisfies the following conditions.
\begin{itemize}
\item
$L_{\infty} = (B, t) + R(\infty)$ is a trivial knot in $S^3$.

\item
$L_s = (B, t) + R(s)$ is a Montesinos link $M(R_1, \ldots, R_k)$. 
\end{itemize}

As in Section~\ref{section:covering},
let $\pi_{\infty}: S^3 \to S^3$ (resp.\
$\pi_s: X_s\to S^3$) be the $2$--fold cover of $S^3$
branched along $L_{\infty}$ (resp.\ $L_s$).
We denote by $K$
the covering knot of the trivializable tangle $(B, t)$,
and by $\gamma_s$ the covering slope of
the $s$--untangle surgery on $L_{\infty}$. 
Then $(K, \gamma_s)$ is a Seifert fibered surgery.
The Montesinos link 
$L_s$ can be deformed into a standard position as in
Figure~\ref{Montesinos}(2).
We define a leading arc of a rational tangle.
Then we show that the preimage of a leading arc becomes a Seifert fiber in $X_s$.

\begin{figure}[h]
\begin{center}
\includegraphics[width=0.7\linewidth]{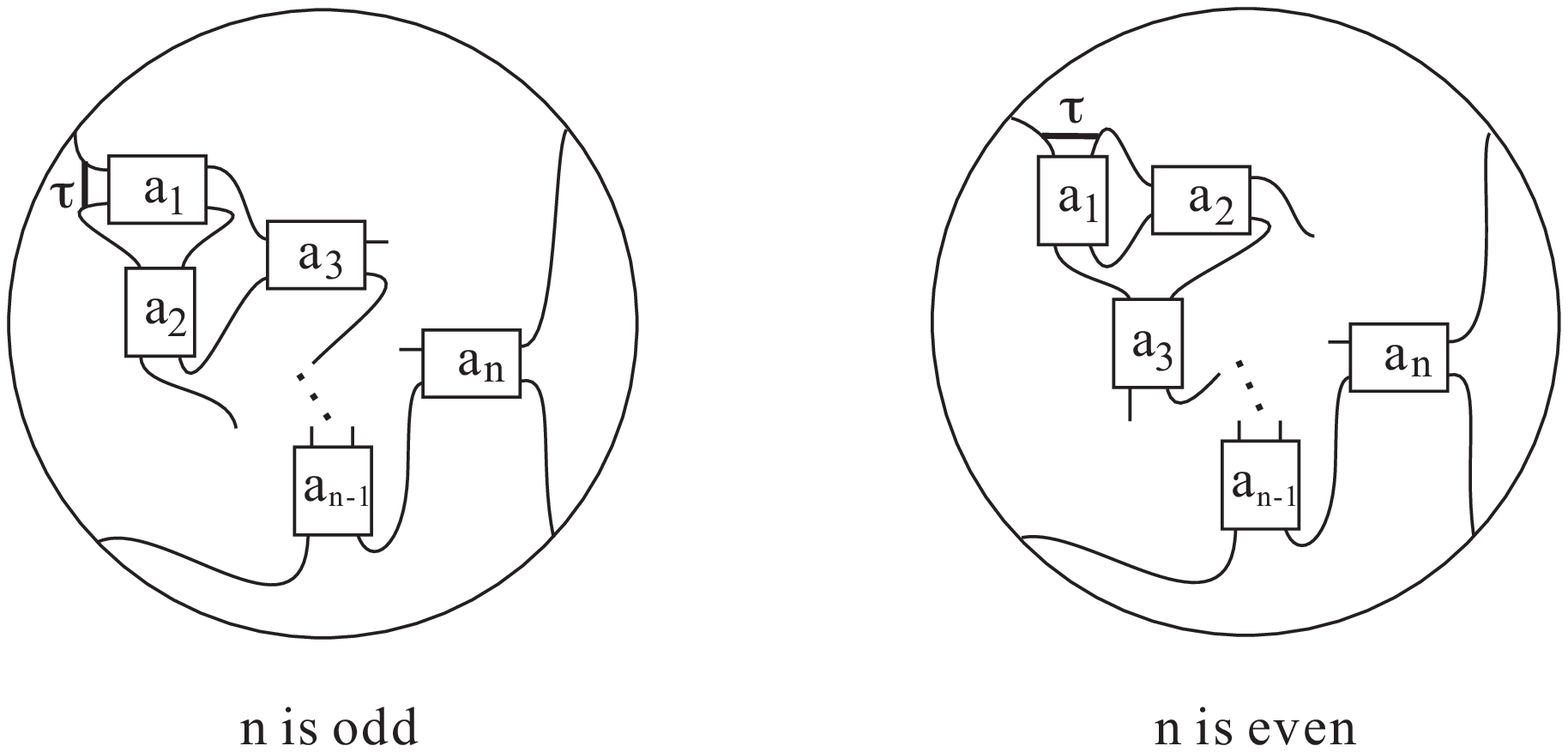}
\caption{Leading arcs in rational tangles.}
\label{fig:leadingarc}
\end{center}
\end{figure}

\begin{definition}[leading arc]
\label{leadingarc}
Let $\tau$ be an arc in a rational tangle
$R = R(a_1, \dots, a_n)$ 
as depicted in Figure~\ref{fig:leadingarc}. 
Then we call $\tau$ a {\it leading arc} of $R$.
\end{definition}

\begin{lemma}
\label{fiber}
Let $\tau$ be a leading arc of a rational tangle
$R_i = R(\frac{p_i}{q_i}) = (B_i, t_i)$ in a standard position of $L_s$.
Then $c = \pi_s^{-1}(\tau)$ is the core of the solid torus
$\pi_s^{-1}(B_i)$ and a fiber of index $|q_i|$
in a Seifert fibration of $X_s$. 
\end{lemma}

\textsc{Proof of Lemma~\ref{fiber}.}
By an ambient isotopy of $B_i$ we can deform $(B_i, t_i)$ and $\tau$ 
to $R(\infty)$ and a spanning arc for $R(\infty)$
as in Figure~\ref{spanningarc}.
Hence $c$ is the core of the solid torus $\pi_s^{-1}(B_i)$.
This implies the desired result.
\QED{Lemma~\ref{fiber}}

\begin{remark}
\label{indexp} 
In Lemma~\ref{fiber}, 
if $|q_i|=1$, 
then $c$ is a regular fiber in $X_s$. 
If $q_i = 0$ i.e.\ $R_i = R(\infty)$, 
then $L_s$ is not a Montesinos link in the usual sense, 
and $X_s$ is a connected sum of lens spaces having
a Seifert fibration with $c$ a degenerate fiber.
\end{remark}

For an arc $\tau$ with $\tau \cap L_{\infty} = \partial \tau$
we perform an untangle surgery along $\tau$ as follows.
First take a regular neighborhood $N(\tau)$ of $\tau$ so that
$T = (N(\tau), N(\tau) \cap L_{\infty})$ is a trivial tangle.
Then, identifying $T$ with the rational tangle $R(\infty)$,
replace $T$ by a rational tangle $R(s)$;
this operation is called $s$--untangle surgery on $L_{\infty}$
along $\tau$.
Then, performing $s$--untangle surgery along $\tau$ downstairs
corresponds to performing Dehn surgery on
the knot $\pi_{\infty}^{-1}( \tau )$ upstairs.

\begin{theorem}[{\cite[Theorem~3.4]{DEMM}}]
\label{seiferter for covering knots}
Let $\tau$ be an arc in $\mathrm{int}B$ such that $\tau \cap t = \partial \tau$.  
Assume that after an isotopy of $\tau \cup L_s$, 
$\tau$ is a leading arc of some $R_i$
in a standard position of $L_s$.  
Assume also that some nontrivial untangle surgery on $L_{\infty}$
along $\tau$ preserves the triviality of $L_{\infty}$.
Then the following hold. 
\begin{enumerate}
\item
The preimage $c = \pi_{\infty}^{-1}(\tau)$ is
a seiferter for $(K, \gamma_s)$.

\item The above untangle surgery along $\tau$ corresponds to twisting along the seiferter $c$ in $S^3$.
\end{enumerate}
\end{theorem}

\begin{remark}
\label{framing}
Assume that $\frac{1}{n_0}$--untangle surgery along $\tau$ preserves
the triviality of $(B, t) + R(\infty)$
for some integer $n_0$ with $|n_0| > 2$. 
Then the preimage of
the latitude of $R(\infty)=(N(\tau), N(\tau)\cap L_{\infty})$ 
is a preferred longitude of $c$ \cite[Remark~3.5]{DEMM}. 
Thus for any integer $n$, 
$\frac{1}{n}$--untangle surgery along $\tau$ corresponds to 
$(-\frac{1}{n})$--surgery on $c$,
i.e.\ $n$--twist along $c$. 
\end{remark}

\subsection{Locating \bm{$(K(l,m,n,p), \gamma_{l,m,n,p})$}
in the Seifert Surgery Network}
\label{Location}

Let $B(l,m,n,p)$, where $m$ or $p$ is 0, be the trivializable tangle
in Figure~\ref{Blmnp},
and $a$ and $b$ the arcs depicted in Figure~\ref{Blmnp}.
Recall that $K(l,m,n,p)$ and $\gamma_{l,m,n,p} =\gamma$ are
the covering knot of $B(l,m,n,p)$ in Figure~\ref{Blmnp}
and the covering slope of $1$--untangle surgery on
$B(l,m,n,p)+R(\infty)$ respectively,
and $B(l,m,n,p) +R(1)$ is a Montesinos link 
(Lemma~\ref{Blmnp+R}).
We denote the $2$--fold cover branched along $B(l,m,n,p)+R(\infty)$ 
(resp.\ $B(l,m,n,p)+R(1)$) by
$\pi_{\infty} : S^3 \to S^3$
(resp.\ $\pi: K(l,m,n,p)(\gamma) \to S^3$).
Set $c_a(l,m,n,p) = \pi_{\infty}^{-1}(a)$ and
$c_b(l,m,n,p) = \pi_{\infty}^{-1}(b)$.
For simplicity, we often write $c_a$, $c_b$ for
$c_a(l,m,n,p)$, $c_b(l,m,n,p)$.

\begin{proposition}
\label{seiferter for Klmnp}
\begin{enumerate}
\item
If $p = 0$, then $c_a(l,m,n,0)$ is
a seiferter for $(K(l, m, n, 0), \gamma)$. 
\item
If $m = 0$, then $c_b(l,0,n,p)$ is 
a seiferter for $(K(l, 0, n, p), \gamma)$. 
\end{enumerate}
\end{proposition}

\begin{figure}[h]
\begin{center}
\includegraphics[width=0.35\linewidth]{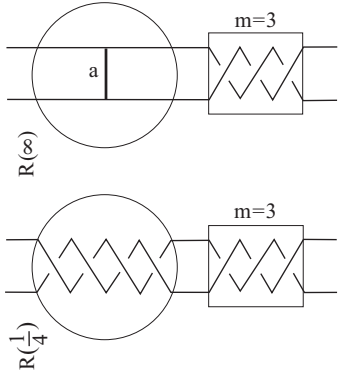}
\caption{$\frac{1}{4}$--untangle surgery along $a$.}
\label{untangle}
\end{center}
\end{figure}

\textsc{Proof of Proposition~\ref{seiferter for Klmnp}.}
$(1)$ Figure~\ref{Blmn0+1Montesinos} gives a standard position of
the Montesinos link $B(l, m, n, 0) + R(1)$ and
shows that the arc $a$ is a leading arc of
$\displaystyle  R_1 = (B_1, t_1)$.
Apply $\frac{1}{m'}$--untangle surgery on $B(l, m, n, 0) + R(\infty)$ along $a$ as in Figure~\ref{untangle}. 
We then obtain $B(l, m-m', n, 0) + R(\infty)$, 
which is a trivial knot by Lemma~\ref{Blmnp+R}(1)(i). 
It follows from Theorem~\ref{seiferter for covering knots}(1) 
that $c_a$ is a seiferter for $(K(l, m, n, 0), \gamma)$.   
In particular, 
by Lemma~\ref{fiber} $c_a$ is the core of $\pi^{-1}(B_1)$
 and an exceptional fiber in $K(l, m, n, 0)(\gamma)$. 

$(2)$ Figure~\ref{Bl0np+1Montesinos} gives a standard position of 
the Montesinos link $B(l, 0, n, p) + R(1)$ and
shows that the arc $b$ is
a leading arc of $\displaystyle R_2 = (B_2, t_2)$.
Note that $\frac{1}{p'}$-untangle surgery on
$B(l, 0, n, p) + R(\infty)$ 
along the arc $b$ as in Figure~\ref{untangle} yields the trivial knot 
$B(l, 0, n, p-p') + R(\infty)$ (Lemma~\ref{Blmnp+R}(2)(i)). 
It follows from Theorem~\ref{seiferter for covering knots}(1) that 
$c_b$ is a seiferter for $(K(l, 0, n, p), \gamma)$; 
in particular, $c_b$ is the core of $\pi^{-1}(B_2)$ and
an exceptional fiber in $K(l, m,0,p)(\gamma)$.  
This establishes Proposition~\ref{seiferter for Klmnp}. 
\QED{Proposition~\ref{seiferter for Klmnp}}

\begin{remark}
\label{remark:c_a c_b}
Figure~\ref{Blmn0+1Montesinos} shows that
the arc $b$ is isotopic to a leading arc of $R_2$
in a standard position of the Montesinos link $B(l,m,n,0)+R(1)$.
This implies that $c_a$ and $c_b$ become
fibers in $K(l,m,n,0)(\gamma_{l,m,n,0})$ simultaneously.
However, if $m\ne 0$, $c_b$ is not necessarily a trivial knot.
If $m=p=0$, then
$c_a$ and $c_b$ are trivial knots, and
$\{c_a, c_b\}$ is a pair of seiferters
for $(K(l,0,n,0), \gamma_{l,0,n,0})$
\end{remark}

Since $\frac{1}{m'}$--untangle surgery on $B(l, m, n, 0) + R(\infty)$ along the arc $a$
preserves the triviality for any $m'$, 
by Remark~\ref{framing}
the $\frac{1}{m'}$--untangle surgery
corresponds to $m'$--twist along the seiferter $c_a$.
Since untangle surgeries along the arc $a$
do not affect the attached tangle $R(\infty)$
in $B(l, m, n, 0) + R(\infty)$,
the image of the covering slope $\gamma_{l, m, n, 0}$
under $m'$--twist along $c_a$
corresponds to 1--untangle surgery on
$B(l, m-m', n, 0) +R(\infty)$.
Thus, $m'$--twist along $c_a$ converts
$(K(l, m, n, 0), \gamma_{l, m, n, 0})$ to
$(K(l, m-m', n, 0), \gamma_{l, m-m', n, 0})$.
Note also that
$m'$--twist along $c_a(l,m,n,0)$ converts
the link $c_a(l,m,n,0) \cup c_b(l,m,n,0)$ to
$c_a(l,m-m',n,0) \cup c_b(l,m-m',n,0)$.
Similar results hold for $p'$--twist of
$(K(l, 0, n, p), \gamma_{l,0,n,p})$ along $c_b$.
We thus have Lemma~\ref{Kl0n0} below.

\begin{lemma}
\label{Kl0n0}
\begin{enumerate}
\item
$m$--twist along $c_a$ converts
$(K(l, m, n, 0), \gamma_{l, m, n, 0})$ to
$(K(l, 0, n, 0), \gamma_{l, 0, n, 0})$. 

\item
$p$--twist along $c_b$ converts 
$(K(l, 0, n, p), \gamma_{l, 0, n, p})$ to 
$(K(l, 0, n, 0), \gamma_{l, 0, n, 0})$. 

\end{enumerate}
\end{lemma}

\begin{figure}[h]
\begin{center}
\includegraphics[width=0.45\linewidth]{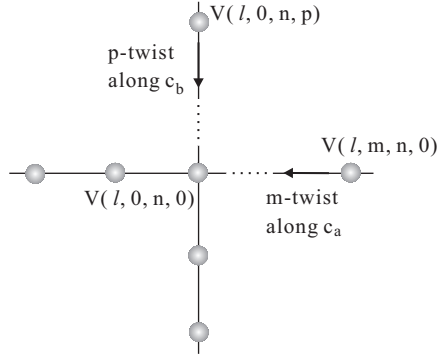}
\caption{$V(l, m, n, p) = (K(l, m, n, p), \gamma_{l,m,n,p})$.}
\label{lattice1}
\end{center}
\end{figure}

Lemma~\ref{Kl0n0}(1) and (2) give the horizontal line and 
the vertical line in Figure~\ref{lattice1}, respectively.  

\begin{lemma}
\label{isotopyB}
\begin{enumerate}
\item
Set $B(l, 1, n-1, 0) =(B, t_1)$ and
$B(l, 0, n, 1) =(B, t_2)$.
Then an ambient isotopy of $B$ fixing $\partial B$ sends
$t_1$ to $t_2$, and the arcs $a$, $b$ to $a$, $b$, respectively.
\item
$(K(l, 1, n-1, 0), \gamma_{l, 1, n-1, 	0}) 
= (K(l, 0, n, 1), \gamma_{l, 0, n, 1})$.
Moreover,
the ordered link $K(l,1,n-1,0) \cup c_a(l,1,n-1,0) \cup c_b(l,1,n-1,0)$
is isotopic to
$K(l,0,n,1) \cup c_a(l,0, n, 1) \cup c_b(l,0, n, 1)$,
so that $\{c_a, c_b\}$ is a pair of seiferters for 
$(K(l, 1, n-1, 0), \gamma)$.

\end{enumerate}
\end{lemma}

\begin{figure}[h]
\begin{center}
\includegraphics[width=0.8\linewidth]{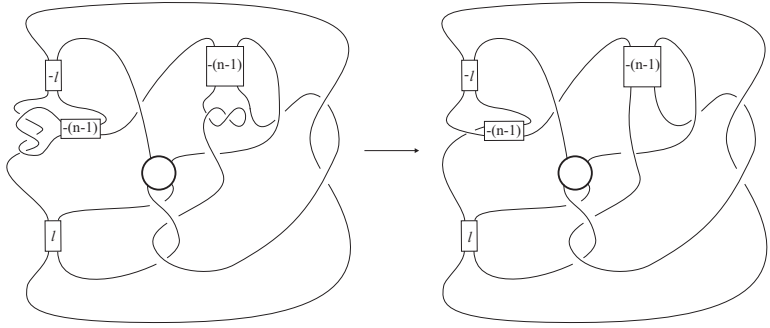}
\caption{An isotopy of $B(l, 1, n-1, 0)$.}
\label{Bl1n0}
\end{center}
\end{figure}

\begin{figure}[h]
\begin{center}
\includegraphics[width=1.0\linewidth]{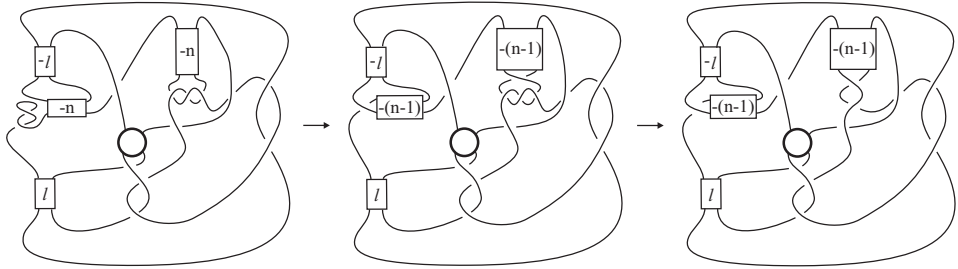}
\caption{An isotopy of $B(l, 0, n, 1)$.}
\label{Bl0n+11}
\end{center}
\end{figure}

\textsc{Proof of Lemma~\ref{isotopyB}.}
$(1)$ We give a pictorial proof.  
The right-most figures in
Figures~\ref{Bl1n0} and \ref{Bl0n+11} depict
the same tangle with arcs $a$, $b$.
Hence, the isotopies in Figures~\ref{Bl1n0} and \ref{Bl0n+11}
imply that an ambient isotopy of $B$ fixing $\partial B$ converts
$B(l, 1, n-1, 0)$ to $B(l, 0, n, 1)$,
and sends the arcs $a$, $b$ to $a$, $b$, respectively.

$(2)$ 
The isotopy in Assertion~(1) extends to
an ambient isotopy of $S^3$ which sends
$B(l, 1, n-1, 0)+R(\infty)$ to $B(l, 0, n, 1)+R(\infty)$,
and the arcs $a$, $b$ to $a$, $b$, respectively.
Hence, there is an ambient isotopy of $S^3$
which sends the ordered link
$K(l, 1, n-1, 0) \cup c_a (l,1,n-1,0) \cup c_b(l,1,n-1,0)$
to $K(l, 0, n, 1) \cup c_a(l,0, n, 1) \cup c_b(l, 0, n, 1)$,
and the covering slope $\gamma_{l, 1, n, 0}$ to
$\gamma_{l, 0, n, 1}$, as claimed.
\QED{Lemma~\ref{isotopyB}}

\begin{figure}[h]
\begin{center}
\includegraphics[width=0.9\linewidth]{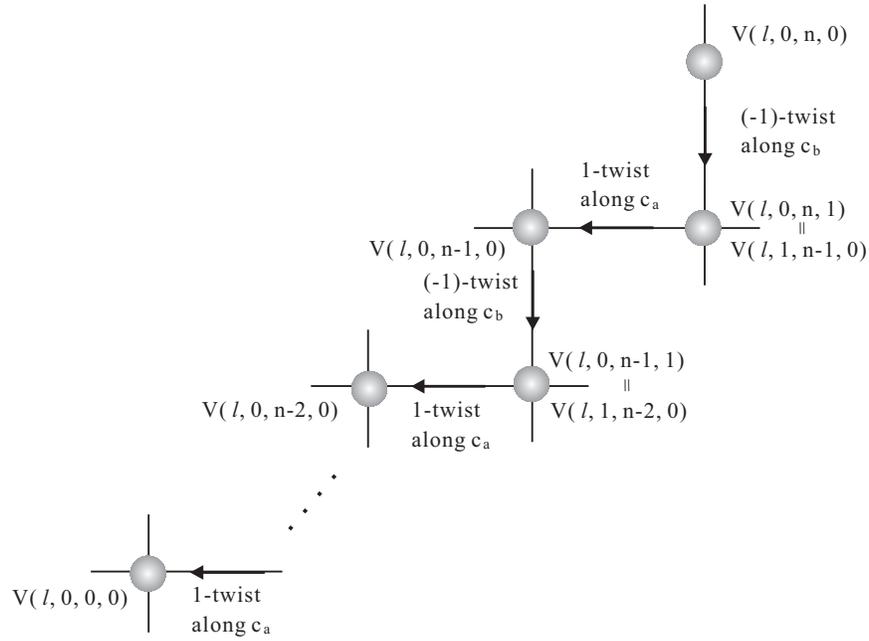}
\caption{$V(l, m, n, p) = (K(l, m, n, p), \gamma_{l,m,n,p})$.}
\label{Kl0n0lattice}
\end{center}
\end{figure}

Applying Lemmas~\ref{Kl0n0}, \ref{isotopyB}(2) repeatedly, 
we find a path 
from $(K(l, 0, n, 0), \gamma_{l,0,n,0})$ to 
$(K(l, 0, 0, 0), \gamma_{l,0,0,0})$
as in Figure~\ref{Kl0n0lattice}. 
Joining this path and the path in Figure~\ref{lattice1}
gives an explicit path from $(K(l,m,n,p), \gamma)$
($m$ or $p$ is $0$) to $(K(l,0,0,0), \gamma)$.

Now we identify $(K(l, 0, 0, 0), \gamma)$
and its seiferters $c_a$, $c_b$.  

\begin{lemma}
\label{trefoil}
$K(l,0,0,0)$ is the trefoil knot $T_{3, 2}$, 
$\gamma_{l,0,0,0} = l +5$,
and $\{ c_a, c_b\}$ is an annular pair of seiferters 
for $(T_{3, 2}, l+5)$ which form
the $(4, 2)$ torus link.
Furthermore, the pair of seiferters $\{c_a, c_b\}$
is the mirror image of
the pair of seiferters $\{c_1^m, s_{-3}\}$
for $(T_{-3, 2}, m)$ given in \cite[Figure~4.2]{DMMtrefoil}
with $m =-l-5$.
\end{lemma}

Combining Figures~\ref{lattice1}, \ref{Kl0n0lattice}, and
Lemma~\ref{trefoil},
we obtain the following result.

\begin{proposition}
\label{location of Klmnp}
The Seifert fibered surgery $(K(l, m, n, p), \gamma_{l,m,n,p})$,
where $m$ or $p$ is $0$,
is obtained from $(T_{3, 2}, l+5)$ by 
a sequence of twists along the pair of seiferters $c_a, c_b$
depicted in Figure~\ref{annularpair_Kl0002_2}
as follows:
alternate $2n$ twists $(-1)$--twist along $c_a$,
$1$--twist along $c_b$,
\ldots, $(-1)$--twist along $c_a$, $1$--twist along $c_b$ $($Figure~\ref{Kl0n0lattice}$)$, 
and finally $(-m)$--twist along $c_a$ or $(-p)$--twist along $c_b$
according as $p=0$ or $m=0$ $($Figure~\ref{lattice1}$)$. 
\end{proposition}

\begin{figure}[h]
\begin{center}
\includegraphics[width=0.48\linewidth]{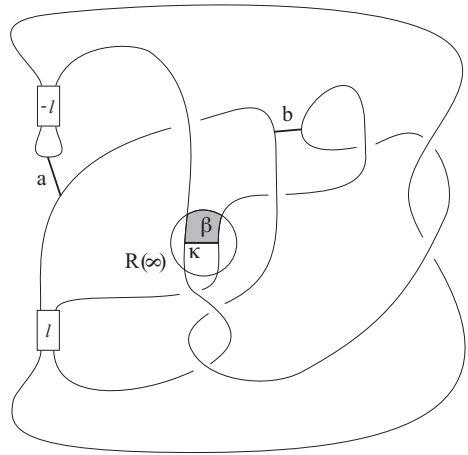}
\caption{$B(l,0,0,0)+R(\infty) \cup a \cup b \cup \kappa$
and band $\beta$.}
\label{Bl000+8}
\end{center}
\end{figure}

\begin{figure}[h]
\begin{center}
\includegraphics[width=1.0\linewidth]{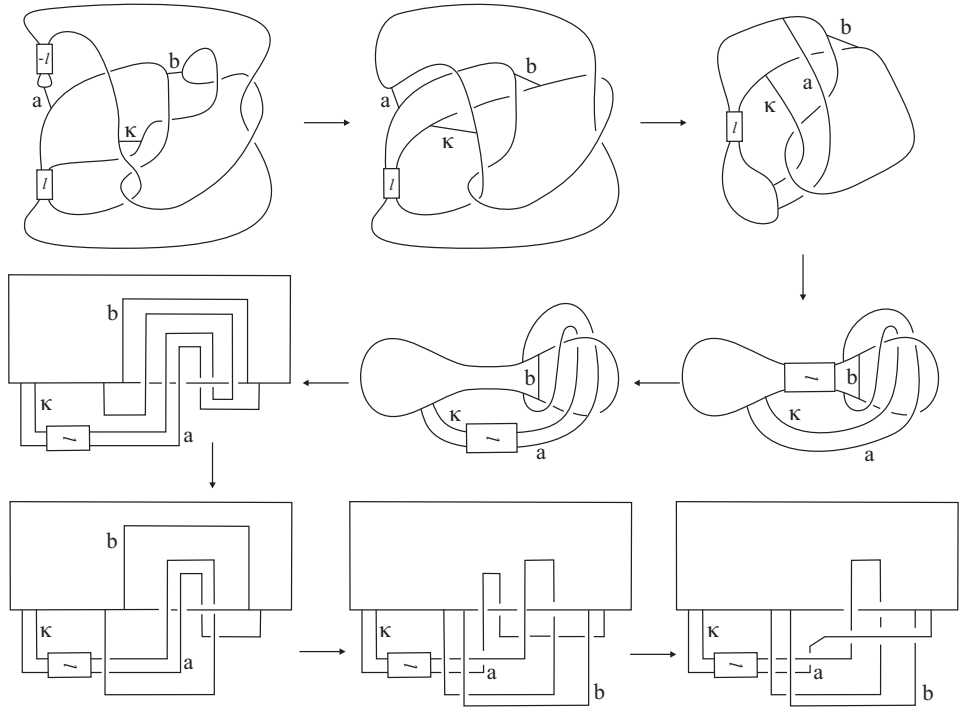}
\caption{An isotopy of $B(l,0,0,0)+R(\infty) \cup a \cup b \cup \kappa$.}
\label{annularpair_Kl000}
\end{center}
\end{figure}

\begin{figure}[h]
\begin{center}
\includegraphics[width=0.24\linewidth]{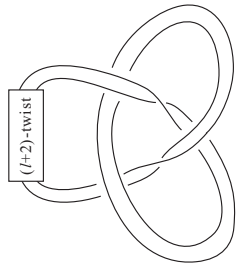}
\caption{Twisted annulus $\pi_{\infty}^{-1}(\beta)$.}
\label{twisted_annulus}
\end{center}
\end{figure}

\begin{figure}[h]
\begin{center}
\includegraphics[width=1.0\linewidth]{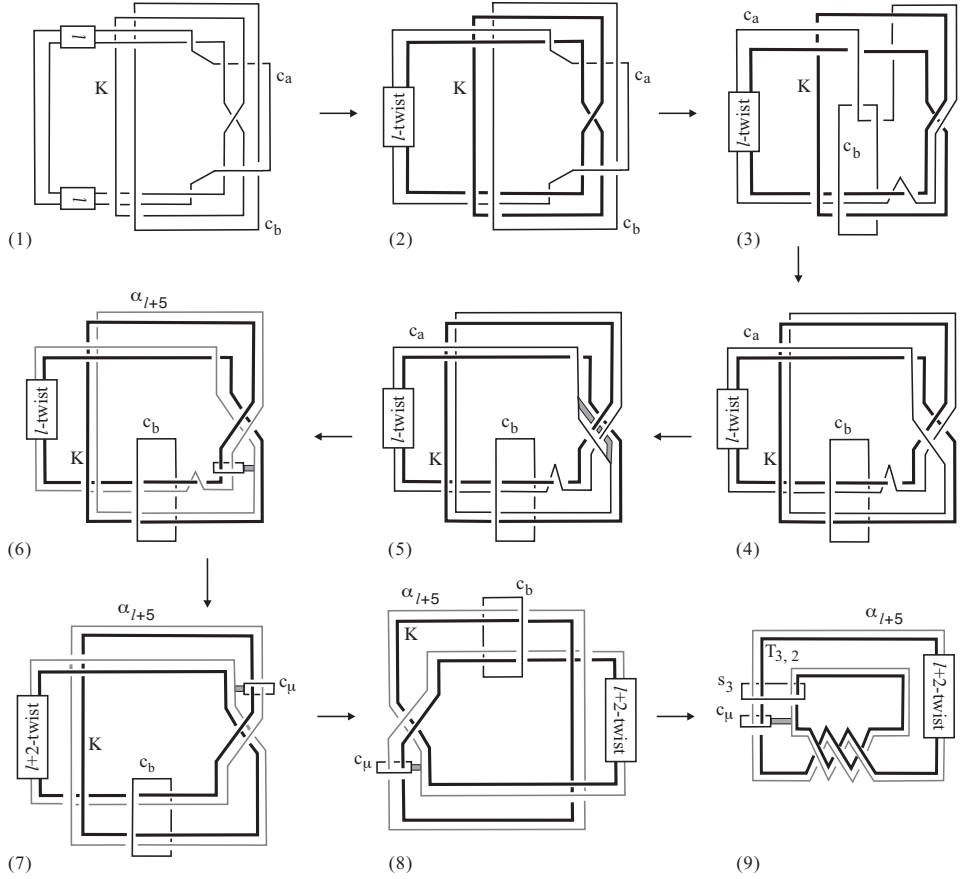}
\caption{$K = K(l,0,0,0)$, $c_a$, and $c_b$.}
\label{annularpair_Kl0002_2}
\end{center}
\end{figure}

\noindent
\textsc{Proof of Lemma~\ref{trefoil}.} 
Figure~\ref{Bl000+8} illustrates the trivial knot 
$B(l, 0, 0, 0) + R(\infty)$, the arcs $a$, $b$,
and the band $\beta$;
$\partial \beta$ is the union
of the spanning arc $\kappa$, 
the latitude of $R(\infty)$,
and two subarcs of $B(l, 0, 0, 0)+R(\infty)$. 
Figure~\ref{annularpair_Kl000} gives an isotopy
of $(B(l, 0, 0, 0) + R(\infty)) \cup \kappa \cup a \cup b$ so that 
$B(l, 0, 0, 0) + R(\infty)$ becomes a standardly embedded circle.
Isotope also the band $\beta$ in the same manner
as in Figure~\ref{annularpair_Kl000}.
Then $(\frac{l}{2} +1)$--twist is added to the band.
Now we consider the $2$--fold branched cover 
$\pi_{\infty}: S^3 \to S^3$ along $B(l, 0, 0, 0) + R(\infty)$.
The preimage $\pi_{\infty}^{-1}( \beta )$ is 
the twisted annulus in Figure~\ref{twisted_annulus}.
Note that the linking number of the parallelly oriented
boundary components of the annulus is 
$2(\frac{l}{2} + 1) +4 = l+6$.
That is, the preimage of the latitude of $R(\infty)$
is a longitude of the solid torus $N(\pi_{\infty}^{-1}(\kappa))$
giving $(l + 6)$--framing.
Then, by Remark~\ref{n--framing}
the covering slope $\gamma_{l,0,0,0}$
of 1--untangle surgery on $B(l,0,0,0)$ equals $l +5$.

The preimages of $\kappa$, $a$, and $b$ become
$K =K(l, 0, 0, 0)$, $c_a$, and $c_b$ in the first figure of 
Figure~\ref{annularpair_Kl0002_2}.
The covering knot $K$ is the trefoil knot $T_{3, 2}$,
and $c_a \cup c_b$ is the $(4, 2)$ torus link
bounding an annulus;
$\{c_a, c_b\}$ is an annular pair of seiferters
for $(T_{3,2}, l+5)$.
We isotope $K\cup c_a\cup c_b$ as in Figure~\ref{annularpair_Kl0002_2}.
Figure~\ref{annularpair_Kl0002_2}(4) shows that
$c_b$ is the exceptional fiber of index $3$
in $S^3 - \mathrm{int}N(T_{3,2})$.
We see from (5) and (6) in Figure~\ref{annularpair_Kl0002_2} that
$c_a$ is a band sum of a knot $c_{\mu}$ in
$S^3 -N(T_{3, 2})$ and
a simple closed curve $\alpha_{l+5}$ in $\partial N(T_{3, 2})$,
where $c_{\mu}$ is parallel to a meridian of $N(T_{3, 2})$
and the slope of $\alpha_{l+5}$ in $\partial N(T_{3, 2})$ is $l+5$.
The last figure of Figure~\ref{annularpair_Kl0002_2} shows that
the mirror image of the ordered link $T_{3, 2} \cup c_a \cup c_b$
is $T_{-3, 2} \cup c_1^m \cup s_{-3}$
depicted in \cite[Figure~4.2]{DMMtrefoil} with $m =-l-5$.
\QED{Lemma~\ref{trefoil}}

As Figure~\ref{Kl0n0lattice} shows,
$K(l,0,n,0)$ is obtained from $K(l,0,0,0)$
by applying a pair of successive twists
$(-1)$--twist along $c_a$ and $1$--twist along $c_b$,
repeatedly $n$ times.
Under $(-1)$--twist along $c_a$ and then $1$--twist along $c_b$
the $(4, 2)$ torus link $c_a \cup c_b$ changes
first to the $(-4, 2)$ torus link and then to the $(4, 2)$ torus link.
We show that applying this sequence of twists
is equivalent to twisting along the annulus cobounded by
$c_a \cup c_b$.
Hence, the Seifert fibered surgery $(K(l, 0, n, 0), \gamma)$
is obtained from $(T_{3, 2}, l +5)$ by twisting along
the annular pair of seiferters $\{c_a, c_b\}$.

\begin{definition}[twist along an annular pair]
\label{def:annulus twist}
Let $c_1, c_2$ be knots in $S^3$
cobounding an annulus $A$,
and give orientations to $c_1, c_2$ so that they are
homologous in $A$.
We call the ordered pair $(c_1, c_2)$ an \textit{annular pair}.
A $p$--\textit{twist along an annular pair} $(c_1, c_2)$ is defined to
be performing $(-\frac{1}{p} +l)$--surgery along $c_1$
and simultaneously $(\frac{1}{p} +l)$--surgery along $c_2$,
where $l =\mathrm{lk}(c_1, c_2)$.
\end{definition}

\begin{lemma}
\label{lem:annulus twist}
Let $c_1\cup c_2$ be the $(4, 2)$ torus link.
Then, performing $(-1)$--twist along $c_1$ and
then $1$--twist along $c_2$ is equivalent to $1$--twist along
the annular pair $(c_1, c_2)$.
\end{lemma}

\textsc{Proof of Lemma~\ref{lem:annulus twist}.} 
The $(4, 2)$ torus link $c_1 \cup c_2$ cobound a unique annulus,
and $l =\mathrm{lk}(c_1, c_2)$ equals 2 when $c_1$ and $c_2$ are
oriented so as to be homologous in the annulus.
Let $\mu$, $\lambda$ be a preferred meridian--longitude pair of $c_2$.
The $(-1)$--twist along $c_1$ changes $m$--framing of $c_2$
to $(m - 2^2)$--framing of $c_2$,
so that $\mu' = \mu$, $\lambda' = \lambda + 4 \mu$ become
 a preferred meridian--longitude pair of
$c_2$ after the twist along $c_1$.
Hence, the surgery slope of
$(-1)$--surgery along $c_2$ (i.e.\ $1$--twist along $c_2$)
after the twist along $c_1$
is $\lambda' - \mu' = \lambda +3\mu$.
Performing $(-1)$--twist along $c_1$ and then $1$--twist along $c_2$
is then performing $1$--surgery along $c_1$
and $3$--surgery along $c_2$ simultaneously.
Since $l =2$,
this shows that the sequence of twists is equivalent to
$1$--twist along the annular pair $(c_1, c_2)$.
\QED{Lemma~\ref{lem:annulus twist}}

\begin{corollary}
\label{cor:location of Klmnp}
The Seifert fibered surgery $(K(l, m, n, p), \gamma_{l,m,n,p})$,
where $m$ or $p$ is $0$,
is obtained from $(T_{3, 2}, l+5)$ by applying
$n$--twist along the annular pair of seiferters $(c_a, c_b)$
depicted in Figure~\ref{annularpair_Kl0002_2}
and then $(-m)$--twist along $c_a$ or $(-p)$--twist along $c_b$
according as $p=0$ or $m=0$.
Regarding the surgery slope $\gamma_{l,m,n,p}$,
$\gamma_{l, 0, n, 0} = 5 +l +n(l^2 +8l +12) +2n^2(l +2)^2$ and
$\gamma_{l, m, n, p} = \gamma_{l, 0, n, 0}
-m(2nl +4n +l +4)^2 -p(2nl +4n +2)^2$.
\end{corollary}

\textsc{Proof of Corollary~\ref{cor:location of Klmnp}.}
The first statement follows from
Proposition~\ref{location of Klmnp} and Lemma~\ref{lem:annulus twist}.
To calculate the surgery slope $\gamma_{l,m,n,p}$
we use results in \cite{DMM1}.
Using Proposition~2.33(2) in \cite{DMM1},
we obtain $\gamma_{l,0,n,0} = \gamma_{l,0,0,0}
+ n( l_1^2 -l_2^2 ) +2n^2(l_1 -l_2)^2$,
where $l_1 = \mathrm{lk}(T_{3, 2}, c_a) = l+4$,
$l_2 = \mathrm{lk}(T_{3, 2}, c_b) =2$
under an adequate orientation of $T_{3, 2}$.
Twisting $(K(l,0,n,0), \gamma_{l,0,n,0})$ $-m$ times along $c_a$
 increases $\gamma_{l,0,n,0}$ by $-m(\mathrm{lk}(K(l,0,n,0), c_a))^2$
\cite[Proposition~2.6]{DMM1}.
Take an annulus cobounded by $c_a \cup c_b$
whose boundary orientation coincides with $c_a$
and the reversed orientation of $c_b$.
Then the annulus is twisted twice, and
intersects $T_{3, 2}$
algebraically $l_1 -l_2 = l+2$ times.
Hence, after $n$--twist along the annular pair $(c_a, c_b)$,
$\mathrm{lk}(K(l,0,0,0), c_a)$ increases by $2n(l+2)$,
so that $\mathrm{lk}(K(l, 0, n, 0), c_a) = l+4 +2n(l +2)$.
This leads to $\gamma_{l,m,n,0}
= \gamma_{l,0,n,0} -m(2nl +4n +l +4)^2$.
Similarly, we obtain the formula of $\gamma_{l,0,n,p}$.
\QED{Corollary~\ref{cor:location of Klmnp}}

\vskip 0.7cm

\noindent
\textbf{Acknowledgements}

\noindent
We would like to thank the referee for careful reading and useful suggestions. 
The first author was partially supported by PAPIIT-UNAM grant IN109811. 
The last author has been partially supported by JSPS Grants-in-Aid for Scientific 
Research (C) (No.21540098), The Ministry of Education, Culture, Sports, Science and Technology, Japan and Joint Research Grant of Institute of Natural Sciences at 
Nihon University for 2013.

\bigskip

\bibliographystyle{amsplain}

\end{document}